\RequirePackage[l2tabu, orthodox]{nag}
\documentclass[a4paper]{amsart}

\usepackage{amscd, amssymb, calc, ifthen, mathrsfs}
\usepackage[matrix,arrow,curve,cmtip]{xy}
\usepackage{amsxtra}
\usepackage{tikz}
\usepackage{microtype}
\usepackage{ellipsis}

\DeclareMathAlphabet{\mathpzc}{OT1}{pzc}{m}{it}

\let \: \colon
\newcommand{\spacesubsec}[1]{\vspace{3mm}\subsec{{#1}}\vspace{2mm}}
\newcommand{\oldmarginpar}[1]{}
\newcommand{\id}{{\mathrm{id}}}

\newcommand{\vbl}{-}
\newcommand{\AngBrak}[1]{{\langle{#1}\rangle}}
\newcommand{\tn}{\otimes}           

\def\isomap{{\,\buildrel \sim\over\rightarrow\,}} 
\def\longisomap{{\,\buildrel \sim\over\longrightarrow\,}} 
\newcommand{\pr}{{\mathrm{pr}}}
\newcommand{\subsec}[1]{{\begin{trivlist}\item\em\large#1\end{trivlist}}}

\newcommand{\longlabelmap}[1]{{\,\buildrel #1\over\longrightarrow\,}}

\newcommand{\longlabelmaps}[2]{{\rightrightarrows}}
\newcommand{\xright}[1]{\xrightarrow{\;#1\;}}

\newcommand{\longmap}{{\,\longrightarrow\,}}

\newcommand{\Hom}{{\mathrm{Hom}}}

\DeclareMathOperator{\Spec}{Spec}

\newcommand{\gr}{{\mathrm{gr}}}

\newcommand{\bF}{{\mathbf{F}}}

\newcommand{\bQ}{{\mathbf{Q}}}

\newcommand{\bZ}{{\mathbf{Z}}}

\newcommand{\bN}{{\mathbf{N}}}

\newcommand{\m}{{\mathfrak{m}}}

\def\longisomap{{\,\buildrel \sim\over\longrightarrow\,}} 

\newcommand{\comment}[1]{}

\renewcommand{\leq}{\leqslant}
\renewcommand{\geq}{\geqslant}

\newcommand{\ord}{{\mathrm{ord}}}
\DeclareMathOperator{\Tor}{Tor}
\DeclareMathOperator{\Sym}{Sym}
\DeclareMathOperator{\colimm}{colim}
\DeclareMathOperator{\limm}{lim}

\newcommand{\smcoprod}{{\,\scriptstyle\amalg\,}}

\newcommand{\eps}{\varepsilon}

\DeclareMathOperator{\im}{im}

\newcommand{\rightlabelxyarrows}[2]{{\ar@<0.7ex>^-{#1}[r]\ar@<-0.7ex>_-{#2}[r]}}

\newcommand{\displaylabelfork}[6]{{	\entrymodifiers={+!!<0pt,\fontdimen22\textfont2>}
	\def\objectstyle{\displaystyle}
\xymatrix{{#1} \ar^-{#2}[r] & {#3} \ar@<0.7ex>^-{#4}[r]\ar@<-0.7ex>_-{#5}[r] & {#6}}}}

\newcommand{\displaylonglabelcofork}[6]{{	\entrymodifiers={+!!<0pt,\fontdimen22\textfont2>}
	\def\objectstyle{\displaystyle}
\xymatrix@C=40pt{{#1} \ar@<0.7ex>^-{#2}[r]\ar@<-0.7ex>_-{#3}[r] & {#4} \ar^-{#5}[r] & {#6}}}}

\newcommand{\displaylabelcofork}[6]{{	\entrymodifiers={+!!<0pt,\fontdimen22\textfont2>}
	\def\objectstyle{\displaystyle}
\xymatrix{{#1} \ar@<0.7ex>^-{#2}[r]\ar@<-0.7ex>_-{#3}[r] & {#4} \ar^-{#5}[r] & {#6}}}}

\newcommand{\displaycofork}[3]{{\displaylabelcofork{#1}{}{}{#2}{}{#3}}}
\newcommand{\displayfork}[3]{{\displaylabelfork{#1}{}{#2}{}{}{#3}}}

\newcommand{\displaylabelrightarrows}[4]{{\entrymodifiers={+!!<0pt,\fontdimen22\textfont2>}
	\def\objectstyle{\displaystyle}
\xymatrix{{#1} \ar@<0.7ex>^-{#2}[r]\ar@<-0.7ex>_-{#3}[r] & {#4}}}}

\newcommand{\displayrightarrows}[2]{{\displaylabelrightarrows{#1}{}{}{#2}}}

\newcommand{\nc}{\!\circ\!}
\newcommand{\setof}[2]{\{#1 \;|\; #2\}}

\newcommand{\bcp}{\odot}
\newcommand{\ptst}{E}
\newcommand{\pind}{\alpha}

\newcommand{\gh}[1]{w_{#1}}
\newcommand{\rgh}[1]{{\bar{w}_{#1}}}

\newcommand{\nset}{{\bN^{(\ptst)}}}
\newcommand{\zset}{{\bZ^{(\ptst)}}}
\newcommand{\EtAlg}{{\mathsf{EtAlg}}}

\newcommand{\wfl}{W^{\mathrm{fl}}}

\newcommand{\ga}{g_1}
\newcommand{\gb}{g_2}
\newcommand{\ha}{h_1}
\newcommand{\hb}{h_2}
\newcommand{\ff}{f}

\newcommand{\PreAct}[2]{\mathrm{PreAct}({#1},{#2})}
\newcommand{\Act}[2]{\mathrm{Act}({#1},{#2})}

\newcommand{\xx}{g}

\newcommand{\cata}{\mathsf{C}}
\newcommand{\catb}{\mathsf{E}}

\newcommand{\pp}[1]{\pi_{#1}}
\newcommand{\inv}{s}

\newcommand{\wpre}{\bar{W}}

\newcommand{\alphapre}{\bar{\alpha}}
\newcommand{\Ring}{\mathsf{Ring}}
\newcommand{\Ringfl}{\mathsf{Ring}^{\mathsf{fl}}}

\newtheoremstyle{mythm}{}{}%
  {\itshape}
  {}
  {\bfseries}
  {}
  { }
  {\thmnumber{#2.\hspace{1.5mm}}\thmname{#1}\thmnote{ #3}.}

\newtheoremstyle{intro}{}{}%
  {\itshape}
  {}
  {\bfseries}
  {}
  { }
  {\thmname{#1}\thmnumber{ #2}\thmnote{ #3}.}

\newtheoremstyle{myrmk}{}{}%
  {}
  {}
  {}
  {}
  { }
  {{\bfseries\thmnumber{#2.\hspace{1.5mm}}}{\itshape\thmname{#1}}\thmnote{ #3}.}

\numberwithin{equation}{subsection}

\theoremstyle{mythm}

\newtheorem{theorem}[subsection]{Theorem}
\newtheorem{prop}[subsection]{Proposition}
\newtheorem{proposition}[subsection]{Proposition}
\newtheorem{lemma}[subsection]{Lemma}

\newtheorem{corollary}[subsection]{Corollary}

\theoremstyle{myrmk}

\theoremstyle{intro}
\newtheorem*{thmintro}{Theorem}

\theoremstyle{plain}
\newtheorem*{prop*}{Proposition}
\newtheorem*{cor*}{Corollary}
\newtheorem*{conj*}{Conjecture}

\theoremstyle{definition}

\makeatletter
\def\@seccntformat#1{\@ifundefined{#1@cntformat}%
{\csname the#1\endcsname\quad}
{\csname #1@cntformat\endcsname}
}
\def\section@cntformat{\thesection.\enspace}
\def\subsection@cntformat{\thesubsection.}
\makeatother

\newcommand\mnote[1]{}

\newcounter{hour}\newcounter{minute}
\newcommand{\printtime}{\setcounter{hour}{\time/60}%
        \setcounter{minute}{\time-\value{hour}*60}%
        \ifthenelse{\value{hour}<10}{0}{}\thehour:%
        \ifthenelse{\value{minute}<10}{0}{}\theminute}

\begin{document}

\title[The basic geometry of Witt vectors, I]{The basic geometry of Witt vectors, I\\The affine case}
\author[J.~Borger]{James Borger}
\address{Australian National University}

\date{\today. \printtime}
\email{james.borger@anu.edu.au}
\thanks{This work was partly supported by Discovery Project DP0773301, 
a grant from the Australian Research Council.}

\begin{abstract}
	We give a concrete description of the category of \'etale algebras over the ring of Witt 
	vectors	of a given finite length with entries in an arbitrary ring.
	We do this not only for the classical $p$-typical and big Witt vector functors
	but also for certain analogues over arbitrary local and global fields. 
	The basic theory of these generalized Witt vectors is developed
	from the point of view of commuting Frobenius lifts and their universal properties,
	which is a new approach even for the classical Witt vectors.
	The larger purpose of this paper is to provide the affine foundations for the algebraic 
	geometry of generalized	Witt schemes and arithmetic jet spaces. So
	the basics here are developed somewhat fully, with an eye toward future applications.
\end{abstract}

\maketitle


\section*{Introduction}

Witt vector functors are certain functors from the category of (commutative) rings to
itself. The most common are the $p$-typical 
Witt vector functors $W$, for each prime number $p$.
Given a ring $A$, one traditionally defines $W(A)$ as a set to be $A^{\bN}$ and then gives it the
unique ring structure which is functorial in $A$ and such that the set maps
	\begin{align*}
	W(A) &\longlabelmap{\gh{}} A^{\bN} \\
	(x_0,x_1,\dots) &\mapsto \AngBrak{x_0,\ x_0^p+px_1,\ x_0^{p^2}+px_1^p+p^2x_2,\ \dots}
	\end{align*}
are ring homomorphisms for all rings $A$, 
where the target has the ring structure with componentwise 
operations. For example, we have
	\begin{align*}
	(x_0,x_1,\dots) + (y_0,y_1,\dots) &= (x_0+y_0,\ 
		x_1+y_1-\sum_{i=1}^{p-1}\frac{1}{p}\binom{p}{i}x_0^iy_0^{p-i},\ \dots) \\
	(x_0,x_1,\dots) \cdot (y_0,y_1,\dots) &= (x_0 y_0,\  
		x_0^p y_1 + x_1 y_0^p + px_1 y_1,\ \dots). \\
	\end{align*}
Observe that the four polynomials in $x_0,y_0,x_1,y_1$ displayed on the right-hand side have
integer coefficients, as they must if they are to define operations on $W(A)$ for all rings $A$.
Conversely, to prove that the desired functorial ring structure on $W$ exists, it is enough to
prove that the polynomials sitting in the higher components have integer 
coefficients too. This is Witt's theorem.

On the other hand, the polynomials at the component of index $n$ depend only on
the variables $x_0,y_0,\dots,x_n,y_n$. This is clear by induction. It follows that the quotient set 
$A^{[0,n]}=\{(x_0,\dots,x_n)\}$ of $W(A)=A^{\bN}$ is
a quotient ring, which we denote by $W_n(A)$. (It is traditionally denoted $W_{n+1}(A)$.
The shift in indexing is preferable for reasons discussed 
in~\ref{subsec:traditional-versus-normalized-indexing}.)

In some cases, the rings $W(A)$ and $W_n(A)$ are isomorphic to familiar rings. For example,
$W(\bZ/p\bZ)$ is isomorphic to the ring $\bZ_p$ of $p$-adic integers, and $W_n(\bZ/p\bZ)$ is
isomorphic to $\bZ/p^{n+1}\bZ$. If $p$ is invertible in $A$, then $\gh{}$ is a bijection and so the
Witt vector rings become product rings: $W_n(A)\cong A^{[0,n]}$ and $W(A)\cong A^{\bN}$. But in 
most cases, $W(A)$ is not a familiar ring.

While this traditional 
approach to Witt vectors is adequate for many purposes, it has two shortcomings. The first
is that it is not clear how we should think about the affine scheme $\Spec W_n(A)$ geometrically. Indeed, I
am not aware of a truly geometric description of $\Spec W_n(A)$ in any nontrivial case in the literature.
If we want to fully incorporate Witt vectors into arithmetic algebraic geometry (and we do), it is
important to have a thorough understanding of their geometry. The main point here and in the
following paper \cite{Borger:BGWV-II} is to set up a framework for that. The geometry in this paper is
however limited to the basic results in the affine case needed for the general treatment in
\cite{Borger:BGWV-II}.

The second shortcoming of the traditional
approach is that it does not explain what the defining purpose of Witt 
vectors is.
The answer, at least for this paper, is that they
control Frobenius lifts---ring endomorphisms which reduce to the Frobenius map modulo $p$.
Here we are following Borger--Wieland \cite{Borger-Wieland:PA}, 12.3--4, which in turn followed
Joyal \cite{Joyal:Lambda}\cite{Joyal:Witt}. Motivated by this perspective, we will define Witt
vector functors relative to primes in any global or local field. This generality includes not only
the $p$-typical functors above but also the so-called big Witt vector functor and less common
variants of the $p$-typical ones due to Drinfeld and to Hazewinkel
(\cite{Drinfeld:symmetric-domains}, Proposition 1.1; \cite{Hazewinkel:book}, (18.6.13)). It also
includes many variants unstudied till now. We will work with these generalized functors
throughout the paper. In fact, this will take no more effort once we establish some basic
reduction results.

\medskip

Let us now discuss the contents in more detail.

Section 1 introduces our generalized Witt vectors. Given a Dedekind 
domain $R$ and a set $\ptst$ of maximal ideals of $R$ with finite residue
fields, we will define a functor $W_{R,\ptst}$ from the category
$\Ring_R$ of $R$-algebras to itself:
	$$
	W_{R,\ptst}\:\Ring_R\to\Ring_R.
	$$
(In fact, we will work with slightly more general $R$ and $\ptst$.) 
We call $W_{R,\ptst}$ the $\ptst$-typical Witt vector functor. When $R=\bZ$ and
$\ptst$ consists of a single maximal ideal $p\bZ$, our functor will agree with the $p$-typical Witt
vector functor above; when $\ptst$ consists of all maximal ideals of $\bZ$, our functor will agree
with the big Witt vector functor. The definition of $W_{R,\ptst}$ is in two steps. First we define
a functor 
	$$
	\wfl_{R,\ptst}\:\Ringfl_R\to\Ringfl_R,
	$$
where $\Ringfl_R$ is the full subcategory of $\Ring_R$ consisting of $R$-algebras which are
$\m$-torsion free for all ideals $\m\in\ptst$. We call such algebras $\ptst$-flat. Then we
define $W_{R,\ptst}$ to be a certain natural extension of $\wfl_{R,\ptst}$ to all of $\Ring_R$.

Let $\nset$ denote the commutative monoid $\bigoplus_{\ptst}\bN$, where $\bN$ is
$\{0,1,\dots\}$ under addition. Given an action of $\nset$ on an $R$-algebra $B$, let $\psi_{\m}$
denote the endomorphism of $B$ given by the $\m$-th element of the standard basis of $\nset$. Let
us say that such an action is a {$\Lambda_{R,\ptst}$-structure} if for each $\m\in\ptst$, the map
$\psi_{\m}$ reduces to the Frobenius endomorphism $x\mapsto x^{[R:\m]}$ on $B/\m B$. Now, for any
$R$-algebra $A$, the monoid $\nset$ acts on $A^{\nset}$ through its translation action on itself in the
exponent. When $A$ is $\ptst$-flat, we define $\wfl_{R,\ptst}(A)$ to be the largest
of the sub-$R$-algebras $B\subseteq A^{\nset}$ having the properties that $B$ is stable under the
action of $\nset$ and that the induced action on $B$ is a $\Lambda_{R,\ptst}$-structure. It is
elementary to check that a maximal such subalgebra $\wfl_{R,\ptst}(A)$ exists.

This definition can be expressed as a universal property. Let
$\Ringfl_{\Lambda_{R,\ptst}}$ denote the following category: the objects are $\ptst$-flat
$R$-algebras equipped with a $\Lambda_{R,\ptst}$-structure, and the morphisms are
$\nset$-equivariant $R$-algebra maps. Then $\wfl_{R,\ptst}$, viewed as a functor
$\Ringfl_R\to\Ringfl_{\Lambda_{R,\ptst}}$, is the right adjoint of the evident forgetful functor.

One then defines $W_{R,\ptst}$ to be the left Kan extension of $\wfl_{R,\ptst}$,
now viewed as a functor $\Ringfl_R\to\Ring_R$. This amounts to the following. It is not hard to
show that the functor $\wfl_{R,\ptst}$ is representable, that is, there exists an $\ptst$-flat
$R$-algebra $\Lambda_{R,\ptst}$ and an isomorphism
$\wfl_{R,\ptst}(\vbl)=\Hom(\Lambda_{R,\ptst},\vbl)$, as set-valued functors. Because 
$\wfl_{R,\ptst}$
takes values in $R$-algebras, $\Lambda_{R,\ptst}$ carries the structure of a co-$R$-algebra object
in $\Ringfl_R$. Because such a structure is described using maps between certain coproducts of
$\Lambda_{R,\ptst}$ with itself, and because $\Ringfl_R$ is a full subcategory of $\Ring_R$ closed
under coproducts, $\Lambda_{R,\ptst}$ continues to be a co-$R$-algebra object when viewed as an
object of $\Ring_R$. Therefore it represents an $R$-algebra-valued functor, and this functor is
what $W_{R,\ptst}$ is defined to be.

Since the definition of $W_{R,\ptst}$ in terms of $\wfl_{R,\ptst}$ is of a purely
category-theoretic nature, one should view the $\ptst$-flat case as the central one.
This is in contrast to the common point of view that the purpose of Witt vector functors is to
lift rings from positive characteristic to characteristic zero.

As in the $\ptst$-flat setting, $W_{R,\ptst}$ is the right adjoint of the forgetful functor
$\Ring_{\Lambda_{R,\ptst}}\to\Ring_R$, but to make sense of this, it is necessary to know the what
a $\Lambda_{R,\ptst}$-structure on a general $R$-algebra is. Unfortunately, it is not easy to
state the definition, and so we will leave it to the body of the paper. In the $\ptst$-flat
setting, it is equivalent to a commuting family of Frobenius lifts indexed by $\ptst$, as above;
but in general, it is a slightly stronger structure that is better behaved. When $R$ is $\bZ$ and
$\ptst$ consists of all maximal ideals of $\bZ$, a $\Lambda_{R,\ptst}$-structure is equivalent to a
$\lambda$-ring structure in the sense of Grothendieck's Riemann--Roch theory, but this
does not admit a simple definition either.

In addition to the right adjoint $W_{R,\ptst}$, the forgetful functor
$\Ring_{\Lambda_{R,\ptst}}\to\Ring_R$ has a left adjoint, which we denote by
$A\mapsto\Lambda_{R,\ptst}\bcp A$.
It has a smaller presence in this paper, but it is very 
important---even in the $p$-typical case, as the work of 
Buium~\cite{Buium:p-jets}\cite{Buium:Arithmetic-diff-equ} makes clear.

Section 2 defines functors $W_{R,\ptst,n}$, which are truncations of $W_{R,\ptst}$ in the same way
that the functors $W_n$ above are truncations of $W$. For any $A\in\Ringfl_R$
and $n\in\nset$, let $\wfl_{R,\ptst,n}(A)$ denote the image of the subring 
$\wfl_{R,\ptst}(A)\subseteq A^{\bN}$ under the canonical projection 
	$$
	A^{\nset}\to A^{[0,n]},
	$$
where 
	$$
	[0,n]=\setof{i\in\nset}{i_{\m}\leq n_{\m}\text{ for all }\m\in\ptst}.
	$$
Then $\wfl_{R,\ptst,n}$ is a functor $\Ringfl_R\to\Ringfl_R$. It is representable by
an $\ptst$-flat $R$-algebra $\Lambda_{R,\ptst,n}$, and we extend it to a functor
	$$
	W_{R,\ptst,n}\:\Ring_R\to\Ring_R
	$$
by taking its left Kan extension, as above.
These truncated functors are related to the original one by the formula
	$$
	W_{R,\ptst}(A) = \lim_n W_{R,\ptst,n}(A).
	$$

Even in the $p$-typical case, this approach to defining the Witt vectors has the advantage over the
traditional one that universal properties are emphasized and the particulars of explicit
constructions are played down. But this comes at a cost. For instance, it is not obvious that
$W_{R,\ptst,n}$ preserves surjectivity of maps. To prove this and other basic facts, it appears
necessary to bring back the Witt components $(x_0,x_1,\dots)$ above, at least in some
form. To define them, the ideals of $\ptst$ must be principal; the purpose of section 3 is to
define them in the minimal case we will need, which is when $\ptst$ consists of a single principal
ideal $\m$. A version of the proof of Witt's theorem then shows there is a unique functorial
bijection $A^{\bN}\to W_{R,\ptst}(A)$ such that when
$A$ is $\ptst$-flat, the composition $A^{\bN}\to W_{R,\ptst}(A) \subseteq A^{\bN}$ satisfies
	$$
	(x_0,x_1,x_2,\dots) \mapsto \AngBrak{x_0,\ x_0^q + \pi x_1,\ 
		x_0^{q^2} + \pi x_1^{q} + \pi^2 x_2,\ \dots}
	$$
where $q=[R:\m]$ and $\pi$ is a fixed generator of $\m$.
We can similarly identify $W_{R,\ptst,n}(A)$ with the quotient $A^{[0,n]}$ consisting of vectors 
$(x_0,\dots,x_n)$. Let me emphasize that the components $(x_0,x_1,\dots)$ depend on the choice of
generator $\pi\in\m$ in a complex, non-multilinear way. But we can use them to define 
Verschiebung operators
	\begin{align*}
		V^j_{\m}\:\m^j\tn_R W_{R,\ptst,n}(A) &\longmap W_{R,\ptst,n+j}(A) \\
		\pi^j \tn (x_0,\dots,x_n) &\mapsto (0,\dots,0,x_0,\dots,x_n),
	\end{align*}
which are independent of the choice of the generator $\pi$.
Making that so is the purpose the tensor factor $\m^j$.

When $\ptst$ consists of a single ideal $\m$ (possibly nonprincipal), section 4 describes 
$W_{R,\ptst,n}$ in terms of the case where $\m$ is principal, which is covered by
section 3. This is done by working Zariski locally on $R$.  
Using the same technique, we will show that the Verschiebung maps as above can be defined when
$\m$ is not assumed to be principal. In fact, there is a unique functorial family of such maps agreeing
with the maps defined above. The image of $V^j_{\m}$ is the kernel of
the canonical projection $W_{R,\ptst,n+j}(A)\to W_{R,\ptst,j}(A)$.

Similarly, section 5 gives a description of $W_{R,\ptst,n}$ when $\ptst$ is general 
in terms of the case where $\ptst$ consists of a single ideal, which is covered by section 4:
if $\m_1,\dots,\m_r$ are the ideals in the support of $n$, then there is a natural isomorphism
	\begin{equation} \label{map:intro-iterate}
		W_{R,\ptst,n} \longisomap W_{R,\m_r,n_{\m_r}}\circ\cdots\circ W_{R,\m_1,n_{\m_1}}.
	\end{equation}
Such a description also holds for $W_{R,\ptst}$, though some care must be taken when $\ptst$ is
infinite. It is also possible to describe the functor $\Lambda_{R,\ptst}\bcp\vbl$, as well as
its truncated variants $\Lambda_{R,\ptst,n}\bcp\vbl$, in terms of the case where
$\ptst$ consists of a single ideal.

Section 6 gives several ring-theoretic facts about $W_{R,\ptst,n}$ which
we will need later. For example, this is where we prove that $W_{R,\ptst,n}$ preserves
surjectivity. Most of the arguments there appear to require the use of Witt components and
the reduction techniques of sections 4 and 5.

\medskip

\medskip

Sections 7--9 prove the main results, which relate Witt vector
functors and \'etale maps. Suppose $\ptst$ consists of a single ideal $\m$. For any ring 
$A$ and any integer $n\geq 1$, we have a diagram
	\begin{equation} \label{diag:intro-gluing}
		\displaylabelfork{W_{R,\ptst,n}(A)}{\alpha_n}{W_{R,\ptst,n-1}(A)\times 
			A}{s\circ\pr_1}{t\circ\pr_2}{A/\m^n A.}
	\end{equation}
When $\m$ is principal, the maps $\alpha_n$, $s$, and $t$ can be defined in terms of the Witt 
components relative to a fixed generator $\pi\in\m$ as follows:
	\begin{align*}
		\alpha_n\:(a_0,\dots,a_n) &\mapsto 
			\big((a_0,\dots,a_{n-1}),a_0^{q^n}+\pi a_1^{q^{n-1}}\dots+\pi^n a_n\big) \\
		s \: (a_0,\dots,a_{n-1}) &\mapsto (a_0^{q^n}+\dots+\pi^{n-1}a_{n-1}^q) \bmod \m^nA \\
		t \: a &\mapsto a\bmod \m^nA.
	\end{align*}
If $A$ is $\m$-torsion free, (\ref{diag:intro-gluing}) is an equalizer diagram.
Figure 1 shows the induced diagram of schemes in the $p$-typical case when $n=1$.

\begin{figure}[t]
\begin{tikzpicture}[scale=0.6]
	\draw (0,0) -- (0,3) -- (4,3) -- (4,0) -- (2,0);
	\draw[dashed] (2,0) -- (0.308,0);
	\draw (0.308,0) -- (0,0);
	\draw (2,0) -- (2,3);

	\draw[dashed] (2,0) -- (3.7,0.616);
	\draw[dashed] (3.7,0.616) -- (3.7,3);
	\draw (3.7,3) -- (3.7,3.616);
	\draw (3.7,3.616) -- (0.308,2.385);
	\draw (0.308,2.385) -- (0.308,-0.616);
	\draw (0.308,-0.616) -- (2,0);

	\draw (8,-2) rectangle (12,1);
	\draw (10,-2) -- (10,1);

	\draw (8,2) rectangle (12,5);
	\draw (10,2) -- (10,5);

	\draw (16,0) -- (16,3);
	
	\path[->] (15.5,2) edge [bend right=15] node[auto,swap] 
		{$\scriptstyle \mathrm{Frobenius}$} (10.5,3.5); 
	\path[->] (15.5,1) edge [bend left=15] node[auto] 
		{$\scriptstyle \mathrm{identity}$} (10.5,-0.5); 
	\path[->] (7,1.5) edge (5,1.5);

	\draw (2.5,-4) node[]{$\Spec W_1(A)$} (9.5,-4) node[]{$\Spec (A \times A)$} 
		(17,-4) node[]{$\Spec A/pA$};
	\path[->] (7,-4) edge (5,-4);
	\path[->] (14.5,-3.85) edge (12.5,-3.85);
	\path[->] (14.5,-4.15) edge (12.5,-4.15);
\end{tikzpicture}
\caption{As a topological space,
$\Spec W_1(A)$ (traditionally denoted $W_2(A)$) is two copies of $\Spec A$ glued along 
$\Spec A/pA$. This is also true as schemes if we assume
that $A$ is $p$-torsion free and we glue transversely and with a Frobenius twist, as indicated.
There is a similar description of $\Spec W_n(A)$ as $\Spec W_{n-1}(A)$ glued with $\Spec A$ 
along $\Spec A/p^nA$.
See the diagram (\ref{diag:intro-gluing}).}
\end{figure}
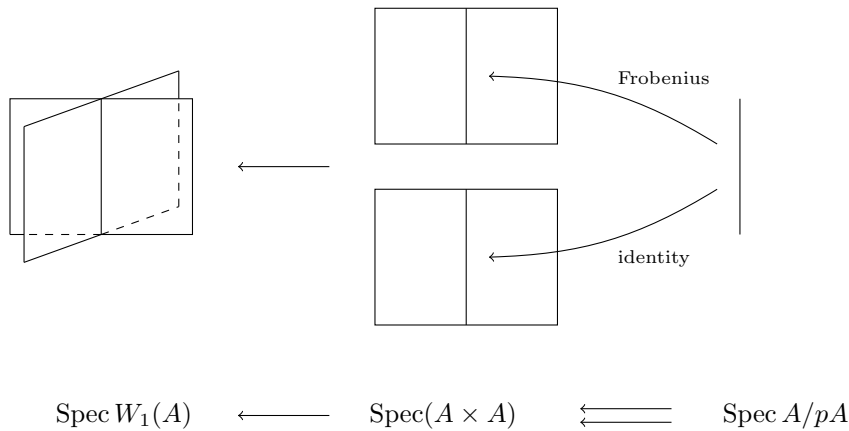

Now let $\mathsf{C}$
denote the following category: an object is a pair $(B,\varphi)$, where $B$ is an \'etale
$(W_{R,\ptst,n-1}(A)\times A)$-algebra and $\varphi$ is an isomorphism of $A/\m^n A$-algebras 
	$$
	A/\m^nA\tn_{t\circ\pr_2}B\longlabelmap{\varphi} A/\m^nA\tn_{s\circ\pr_1}B
	$$
and where a morphism $(B_1,\varphi_1)\to(B_2,\varphi_2)$ 
is a $(W_{R,\ptst,n-1}(A)\times A)$-algebra map $f\:B_1\to B_2$ such that 
	$$
	\varphi_2\circ(A/\m^n A\tn_{t\circ\pr_2}f) 
		= (A/\m^n A\tn_{s\circ\pr_1} f)\circ\varphi_1.
	$$
In other words, $\mathsf{C}$ is the category of algebras equipped with gluing data relative to the
diagram (\ref{diag:intro-gluing}), or equivalently $\mathsf{C}$ is the (weak) fiber product of the
category of \'etale $W_{n-1}(A)$-algebras and the category \'etale $A$-algebras over the category
of \'etale $A/\m^n A$-algebras via the evident functors.

\begin{thmintro}[A]\label{thm:A}
The base-change functor from the category of \'etale $W_{R,\ptst,n}(A)$-algebras to $\mathsf{C}$ is
an equivalence. If $A$ is $\m$-torsion free, then a quasi-inverse is given
by sending $(B,\varphi)$ to the equalizer of the two maps
	$$
	\xymatrix@C=50pt{
	B \rightlabelxyarrows{1\tn\id_B}{\varphi\circ(1\tn\id_B)}
		& A/\m^nA\tn_{s\circ\pr_1}B.
	}
	$$
\end{thmintro}

The first statement can be expressed succinctly in geometric terms; it says that the map $\alpha_n$
satisfies effective descent for \'etale algebras and that descent data is equivalent to gluing data
with respect to the diagram (\ref{diag:intro-gluing}). Using theorem A and induction on $n$, it
is in principle possible to reduce questions about \'etale $W_n(A)$-algebras to questions about
\'etale $A$-algebras. This is still true when $\ptst$ consists of more than one ideal, but now by
(\ref{map:intro-iterate}) and induction on $r$.

Section 9 generalizes van der Kallen's theorem \cite{van-der-Kallen:Descent}, (2.4), 
to arbitrary $R$ and $\ptst$:

\begin{thmintro}[B]\label{thm:B}
	Let $f\:A\to B$ be an \'etale morphism of $R$-algebras. Then the map
	$W_{R,\ptst,n}(f)\: W_{R,\ptst,n}(A)\to W_{R,\ptst,n}(B)$ of $R$-algebras is \'etale.
\end{thmintro}

This theorem is fundamental in extending Witt constructions beyond affine schemes and will be used
often in~\cite{Borger:BGWV-II}. Van der Kallen's argument, which has an infinitesimal flavor, could
be extended to our setting with only minor modifications. Instead we deduce
theorem B from theorem A, and so our argument has a globally geometric flavor. (Note
that, until recently, van der Kallen's paper~\cite{van-der-Kallen:Descent} had escaped the notice
of many workers in de Rham--Witt theory, to whom theorem B was unknown even for the $p$-typical
Witt vectors.)

Last, I would like to thank Amnon Neeman for helpful discussions on some technical points and
Lance Gurney for some comments on earlier versions of this paper.

\tableofcontents

\setcounter{section}{0}

\section{Generalized Witt vectors and $\Lambda$-rings}
\label{sec:generalized-Witt-vectors}

The purpose of this section is to define our generalized Witt vectors and
$\Lambda$-rings. It is largely an expansion in more concrete terms of
Borger--Wieland~\cite{Borger-Wieland:PA}, or rather of the parts about Witt vectors and
$\Lambda$-rings. The approach here will allow us to avoid much of the abstract language of
operations on rings, as first introduced in Tall--Wraith~\cite{Tall-Wraith}.

For the traditional way of defining $\Lambda$-rings and Witt vectors, see \S1 of chapter IX of
Bourbaki~\cite{Bourbaki:CommAlg} and especially the exercises for that section. One can also see
Witt's original paper~\cite{Witt:Vectors} on the $p$-typical Witt vectors (reprinted
in~\cite{Witt:CollectedPapers}) and his notes on the big Witt
vectors~\cite{Witt:CollectedPapers}, pp.\ 157--163.

\subsection{} \emph{Supramaximal ideals.}
Let us say that an ideal $\m$ of a ring $R$ is supramaximal if either
	\begin{enumerate}
		\item $R/\m$ is a finite field, $R_{\m}$ is 
			a discrete valuation ring, and $\m$ is finitely presented as an $R$-module, or
		\item $\m$ is the unit ideal.
	\end{enumerate}
By far the most important example is a maximal ideal with finite residue field in a Dedekind
domain. (In fact, all phenomena in this paper occur already over $R=\bZ$, and this
case covers the classical Witt vectors and $\lambda$-rings.)
The reason we allow the unit ideal is only so that a supramaximal ideal remains
supramaximal after any localization.

Note that a supramaximal ideal $\m$ is invertible as an $R$-module.
Indeed, locally at $\m$ it is the maximal ideal of a discrete valuation ring,
and away from $\m$ it is the unit ideal.

\subsection{} {\em General notation.}
\label{subsec:affine-general-notation}
Fix a ring $R$ and a family $(\m_\pind)_{\pind\in\ptst}$ of pairwise coprime 
supramaximal ideals of $R$ indexed by a set $\ptst$.
Note that because the unit ideal is coprime to itself, it can
be repeated any number of times; otherwise the ideals $\m_{\pind}$ are distinct.
For each $\pind\in\ptst$, let $q_\pind$ denote the cardinality of $R/\m_\pind$.
We will often abusively speak of $\m_\pind$ rather than $\pind$ as being an element of $\ptst$,
especially when $\m_{\pind}$ is maximal, in which case it comes from a unique $\pind\in\ptst$.

Let $R[1/\ptst]$ denote the $R$-algebra whose spectrum is the complement of $\ptst$ in 
$\Spec R$. It is the universal $R$-algebra in which every $\m_\pind$ becomes the 
unit ideal.  It also has the more concrete description 
	$$
	R[1/\ptst]=\bigotimes_{\pind\in\ptst} R[1/\m_{\pind}],
	$$
where the tensor product is over $R$ and $R[1/\m_{\pind}]$ is defined to be
the coequalizer of the maps 
	$$
	\displayrightarrows{\Sym(R)}{\Sym(\m_\pind^{-1})}
	$$
of symmetric algebras, where $\m_\pind^{-1}$ is the dual of $\m_{\pind}$,
one of the maps is $\Sym$ applied to the canonical map $R\to \m_\alpha^{-1}$, and the
other is the map induced by the $R$-module map $R\to\Sym(\m_\pind^{-1})$
that sends $1\in R$ to the element $1\in\Sym(\m_{\pind}^{-1})$ in degree zero.  

Finally, we write $\bN$ for the monoid $\{0,1,2,\dots\}$ under addition and write
$\Ring_R$ for the category of $R$-algebras.

\subsection{} {\em $\ptst$-flat $R$-modules.}
Let us say that an $R$-module $M$ is {\em $\ptst$-flat} if for all
maximal ideals $\m$ in $\ptst$, the following equivalent conditions are satisfied:
	\begin{enumerate}
		\item $R_\m\tn_R M$ is a flat $R_\m$-module,
		\item the map $\m\tn_R M\to M$ is injective.
	\end{enumerate}
The equivalence of these two can be seen as follows.  Condition (b) is equivalent
to the statement $\Tor_1^R(R/\m,M)=0$, which is equivalent to 
	$$
	\Tor_1^{R_\m}(R/\m,R_\m\tn_R M)=0.
	$$
Since $R_\m$ is a discrete valuation ring, this
is equivalent to the $R_\m$-module $R_m\tn_R M$ being torsion free and hence flat.

We say an $R$-algebra is $\ptst$-flat if its underlying $R$-module is. Let $\Ringfl_R$ denote the
full subcategory of $\Ring_R$ consisting of the $\ptst$-flat $R$-algebras.

\begin{proposition}\label{pro:flat=torsion-free}
	Any product of $\ptst$-flat $R$-modules is $\ptst$-flat, and any sub-$R$-module
	of an $\ptst$-flat $R$-module is $\ptst$-flat.
\end{proposition}
\begin{proof}
	We will use condition (b) above.
	Let $(M_i)_{i\in I}$ be a family of $\ptst$-flat $R$-modules.
	We want to show that for each maximal ideal $\m$ in $\ptst$, the
	composition
		\[
		\xymatrix{
		\m\tn\prod_{i} M_i \ar[r] &
			\prod_{i}\m\tn M_i \ar[r] &
			\prod_{i} M_i
		}
		\]
	is injective. Because each $M_i$ is $\ptst$-flat, the right-hand map
	is injective, and so it is enough to show the left-hand map is injective.
	
	Since $\m$ is assumed to be finitely presented as an $R$-module, 
	we can express it as a cokernel of a map $N'\to N$ of finite free
	$R$-modules.  Then we have the following diagram with exact rows
		\[
		\xymatrix{
		N'\tn_R\prod_{i} M_i \ar[r]\ar^\sim[d] &
			N\tn_R\prod_{i}M_i \ar[r]\ar^\sim[d] &
			\m\tn_R \prod_{i} M_i \ar[r]\ar[d] &
			0 \\
		\prod_{i}N'\tn_R M_i \ar[r] &
			\prod_{i} N\tn_R M_i \ar[r] &
			\prod_{i} \m\tn_R M_i \ar[r] &
			0. 
		}
		\]
	The left two vertical maps are isomorphisms because $N'$ and $N$ are finite free.
	Therefore the rightmost vertical map is an injection (and even an isomorphism).
	
	Now suppose $M'$ is a sub-$R$-module of an $\ptst$-flat $R$-module $M$.
	Since $\m$ is an invertible $R$-module, $\m\tn_R M'$ maps injectively to
	$\m\tn_R M$.  Since $M$ is $\ptst$-flat, $\m\tn_R M'$ further maps injectively
	to $M$, and hence to $M'$.
\end{proof}

\subsection{} {\em $\Psi$-rings.}
\label{subsec:Psi-rings}
Let $A$ be an $R$-algebra.  Let us define a $\Psi_{R,\ptst}$-action, or 
a $\Psi_{R,\ptst}$-ring structure, on $A$ to be a commuting family of $R$-algebra
endomorphisms $\psi_\pind$ indexed by $\pind\in\ptst$.  This is the same
as an action of the monoid $\nset=\bigoplus_E\bN$ on $A$.  For any element $n\in\nset$,
we will also write $\psi_n$ for the endomorphism of $A$ induced by $n$.
A morphism of $\Psi_{R,\ptst}$-rings is defined to be an $\nset$-equivariant morphism
of rings.

The free $\Psi_{R,\ptst}$-ring on one generator $e$ is $\Psi_{R,\ptst}=R[e]^{\tn_R\nset}$,
where $\nset$ acts on $\Psi_{R,\ptst}$ through its action on itself in the exponent.
In particular, $\Psi_{R,\ptst}$ is freely generated as an $R$-algebra by the
elements $\psi_n(e)$, where $n\in\nset$. 
Then it is natural to write $\psi_n=\psi_n(e)\in\Psi_{R,\ptst}$ and
$\psi_{\pind}=\psi_{b_{\pind}}\in\Psi_{R,\ptst}$, where $b_\pind\in\nset$ denotes the $\pind$-th
standard basis vector, and $e=\psi_0\in\Psi_{R,\ptst}$ for the identity operator.

For any $\Psi_{R,\ptst}$-ring $A$, there is a unique set map
	\begin{equation} \label{eq:Psi-action}
		\Psi_{R,\ptst}\times A \longlabelmap{\circ} A
	\end{equation}
with the property that for all $\pind\in\ptst$, $r\in R$, $f_1,f_2\in\Psi_{R,\ptst}$, $a\in A$
we have
 	\begin{equation}
	\psi_{\pind}\circ a = \psi_{\pind}(a)
	\end{equation}
and
	\begin{equation}
\quad r\circ a = r, \quad
(f_1+f_2)\circ a = (f_1\circ a) + (f_2\circ a), \quad (f_1f_2)\circ a = (f_1\circ a)(f_2\circ a).
	\end{equation}
Taking $A=\Psi_{R,\ptst}$, we get a binary operation $\circ$ on $\Psi_{R,\ptst}$
called \emph{composition} or \emph{plethysm}. One can check that
this makes $\Psi_{R,\ptst}$ a monoid (noncommutative unless $R=0$) with identity $e$ and that 
(\ref{eq:Psi-action}) is a monoid action. 

In the language of plethystic algebra~\cite{Borger-Wieland:PA}, we can interpret $\Psi_{R,\ptst}$
as the free $R$-plethory $R\AngBrak{\psi_\pind|\pind\in\ptst}$ on the $R$-algebra endomorphisms
$\psi_{\pind}$. Then a $\Psi_{R,\ptst}$-action in the sense above is the same as a
$\Psi_{R,\ptst}$-action in the sense of abstract plethystic algebra. In particular,
$\Psi_{R,\ptst}$ can be viewed as the ring of natural unary operations on $\Psi_{R,\ptst}$-rings,
and the composition operation $\circ$ above agrees 
with the usual composition of unary operations. (Compare with~\ref{subsec:free-lambda-rings} below)

\subsection{} \emph{$\ptst$-flat $\Lambda$-rings.}
\label{subsec:E-flat-lambda-rings}
Let $A$ be an $R$-algebra which is $\ptst$-flat.
Define a $\Lambda_{R,\ptst}$-action, or a $\Lambda_{R,\ptst}$-ring structure,
on $A$ to be a $\Psi_{R,\ptst}$-action with the following 
\emph{Frobenius lift} property: for all $\pind\in\ptst$, 
the endomorphism $\id\tn\psi_{\pind}$ of $R/\m_\pind\tn_R A$ agrees with the Frobenius
map $x\mapsto x^{q_\pind}$.  A morphism of $\ptst$-flat $\Lambda_{R,\ptst}$-rings is 
simply defined to be a morphism of the underlying $\Psi_{R,\ptst}$-rings.
Let us denote this category by $\Ringfl_{\Lambda_{R,\ptst}}$.

\subsection{} \emph{The ghost ring.}
Since an action of $\Psi_{R,\ptst}$ on an $R$-algebra $A$ is the same as an action (in the category
of $R$-algebras) of the monoid $\nset$, the forgetful functor from the category
of $\Psi_{R,\ptst}$-rings to that of $R$-algebras has a right adjoint given by
	$$
	A \mapsto \prod_{\nset}A = A^{\nset},
	$$
where $\nset$ acts on $A^{\nset}$ through its action on itself in the exponent.
(This is a general fact about monoid actions in any category
with products.) For $a\in A^{\nset}$ and $n,n'\in\nset$, the $n$-th component of
$\psi_{n'}(a)$ is the $(n+n')$-th component of $a$.

One might call $A^{\nset}$ the
cofree $\Psi_{R,\ptst}$-ring on the $R$-algebra $A$.
It has traditionally been called the ring of {\em ghost components} or {\em ghost
vectors}. By~\ref{pro:flat=torsion-free}, it is $\ptst$-flat if $A$ is.

When $|E|=1$, there is the possibility of confusing the ghost ring $A^{\bN}$, which has the product
ring structure, with the usual ring $A^{\bN}$ of Witt components (see
\ref{subsec:Witt-components}), which has an exotic ring structure.
To prevent this, we will use angle brackets
$\AngBrak{a_0,a_1,\dots}$ for elements of the ghost ring.

\subsection{} {\em Witt vectors of $\ptst$-flat rings.}
\label{subsec:W-for-flat-rings}
Let us now construct the functor $\wfl_{R,\ptst}$.  We 
will show in~\ref{pro:wfl-univ-property} that it is the right adjoint of the forgetful functor
from the category of $\ptst$-flat $\Lambda_{R,\ptst}$-rings to that of $\ptst$-flat $R$-algebras.
(Further, the flatness will be removed in~\ref{subsec:W-define-general}.)

Let $A$ be a $\ptst$-flat $R$-algebra.  
Let $U_0(A)$ denote the cofree $\Psi_{R,\ptst}$-ring $A^{\nset}$.  For any $i\geq 0$, let
	\[
	U_{i+1}(A) = \setof{b\in U_{i}(A)}
			{\psi_\pind(b)-b^{q_\pind}\in \m_\pind U_{i}(A)\text{ for all } \pind\in\ptst}.
	\]
This is a sub-$R$-algebra of $A^{\nset}$.  Indeed, it is the intersection over $\pind\in\ptst$ 
of the equalizers of pairs of $R$-algebra maps
		\[
		\displayrightarrows{U_i(A)}{R/\m_\pind\tn_R U_i(A)}
		\]
given by $x\mapsto 1\tn \psi_\pind(x)$ and 	by $x\mapsto (1\tn x)^{q_\pind}$.

Now define
	\begin{equation} \label{def:wflat}
		\wfl_{R,\ptst}(A) = \bigcap_{i\geq 0} U_{i}(A).
	\end{equation}
This is the ring of {\em $\ptst$-typical Witt vectors} with entries in $A$. 
It is a sub-$R$-algebra of $A^{\nset}$.
Observe that $\wfl_{R,\ptst}(A)=A^{\nset}$ if $A$ is an $R[1/\ptst]$-algebra.

\begin{proposition}\label{pro:wfl-univ-property}
	\begin{enumerate}
		\item $\wfl_{R,\ptst}(A)$ is a sub-$\Psi_{R,\ptst}$-ring of $A^{\nset}$.
		\item This $\Psi_{R,\ptst}$-ring structure on $\wfl_{R,\ptst}(A)$
			is a $\Lambda_{R,\ptst}$-ring structure.
		\item The induced functor $A\mapsto \wfl_{R,\ptst}(A)$ from $\ptst$-flat 
			$R$-algebras to $\ptst$-flat $\Lambda_{R,\ptst}$-rings is the right adjoint of 
			the forgetful functor.
	\end{enumerate}
\end{proposition}
\begin{proof}
	(a): Let us first show by induction that each $U_i(A)$ is a sub-$\Psi_{R,\ptst}$-ring of 
	$A^{\nset}$.  For $i=0$, we have $U_0(A)=A^{\nset}$, and so it is clear.
	For $i\geq 1$, we use the description of
	$U_{i+1}(A)$ as the intersection 
	of the equalizers of the pairs of ring maps
		\[
		\displayrightarrows{U_i(A)}{R/\m_\pind\tn_R U_i(A)}
		\]
	given in~\ref{subsec:W-for-flat-rings}.
	Observe that both these ring maps become
	$\Psi_{R,\ptst}$-ring maps if we give
	$R/\m_\pind\tn_R U_i(A)$ a $\Psi_{R,\ptst}$-action by defining 
	$\psi_\beta\:a\tn x\mapsto a\tn \psi_\beta(x)$, for all $\beta\in\ptst$.
	Since limits of $\Psi_{R,\ptst}$-rings exist
	and their underlying rings 
	agree with the limits taken in the category of rings, $U_{i+1}(A)$ is a
	sub-$\Psi_{R,\ptst}$-ring of $A^{\nset}$.
	Therefore $\wfl_{R,\ptst}(A)$, the intersection of the $U_i(A)$, is
	also a sub-$\Psi_{R,\ptst}$-ring of $A^{\nset}$.

	(b): It is enough to verify
		$$
		\psi_\pind(x)-x^{q_\pind} \in  
			\m_\pind \wfl_{R,\ptst}(A) = \m_\pind \bigcap_{i\geq 0} U_i(A),
		$$
	for all $\pind\in\ptst$ and $x\in\wfl_{R,\ptst}(A)$.
	For any $i\geq 0$, we know
		$$
		\psi_\pind(x)-x^{q_\pind} \in \m_\pind U_i(A),
		$$
	because $x\in \wfl_{R,\ptst}(A)\subseteq U_{i+1}(A)$. Therefore we know
		$$
		\psi_\pind(x)-x^{q_\pind}\in \bigcap_{i\geq 0} \m_\pind U_i(A).
		$$
	So, it is enough to show
		\begin{equation} \label{eq:internal-nana}
			\m_\pind \bigcap_{i\geq 0} U_i(A) = \bigcap_{i\geq 0} \m_\pind U_i(A). 
		\end{equation}
	Since $\m_\pind$ is finitely generated, it is
	a quotient of a finite free $R$-module $N$. Consider the induced diagram
		$$
		\xymatrix{
		\m_\pind\tn_R \lim_i U_i(A) \ar^{h}[r]
			& \lim_i \m_\pind\tn_R U_i(A)  \\
		N\tn_R \lim_i U_i(A) \ar^{f}[r] \ar[u]
			& \lim_i N\tn_R U_i(A) \ar^{g}[u].
		}
		$$
	Since $N$ is finite free, $f$ is an isomorphism; since $\m_{\pind}$ is projective,
	the map $N\to \m_\pind$ has a section and hence so does $g$. Therefore $g\circ f$ is
	surjective and hence so is $h$, which implies (\ref{eq:internal-nana}).
	
	(c): Let $A$ be an $\ptst$-flat $R$-algebra, 
		let $B$ be an $\ptst$-flat $\Lambda_{R,\ptst}$-ring, and let
		$\bar{\gamma}\:B\to A$ be an $R$-algebra map.  By the cofree property of $A^{\nset}$,
		there is a unique $\Psi_{R,\ptst}$-ring map $\gamma\:B\to A^{\nset}$ lifting 
		$\bar{\gamma}$. We now only
		need to show that the image of $\gamma$ is contained in $\wfl_{R,\ptst}(A)$.
		By induction, it is enough to show that if $\im(\gamma)\subseteq U_i(A)$, then
		$\im(\gamma)\subseteq U_{i+1}(A)$.
		
		Let $b$ be an element of $B$.  Then 
		for each $\pind\in\ptst$, we have
			\[
					\psi_\pind\big(\gamma(b)\big)-\gamma(b)^{q_\pind} 
						= \gamma\big(\psi_\pind(b)-b^{q_\pind}\big) \in \gamma(\m_\pind B)
						\subseteq \m_\pind\im(\gamma) \subseteq \m_\pind U_i(A).
			\]
		Therefore, by definition of $U_{i+1}(A)$, the element $\gamma(b)$ lies in $U_{i+1}(A)$.
\end{proof}

\subsection{} \emph{Exercises.}
\label{subsec:W(Z)-exercise}
Let $R=\bZ$. If $\ptst$ consists of the single ideal $p\bZ$,
then $\wfl(\bZ)$ agrees with the subring
of the ghost ring $\bZ^{\bN}$ consisting of vectors $a=\AngBrak{a_0,a_1,\dots}$ that satisfy 
	$$
	a_{n}\equiv a_{n+1}\bmod p^{n+1}
	$$ 
for all $n\geq 0$.  In particular, the elements are $p$-adic Cauchy sequences
and the rule $a\mapsto \lim_{n\to\infty}a_n$
defines a surjective ring map $\wfl(\bZ)\to\bZ_p$.  

We can go a step further with $\wfl(\bZ_p)$.
Let $I$ denote the ideal $p\bZ_p\times p^2\bZ_p\times\cdots$ in $\bZ_p^{\bN}$.
Then $\wfl(\bZ_p)$ is isomorphic to the ring $\bZ_p \oplus I$
with multiplication defined by the formula $(x,y)(x',y')=(xx',xy'+yx'+yy')$.

Now suppose that $\ptst$ consists of all the maximal ideals of $\bZ$,
and identify $\nset$ with the set of positive integers, by unique factorization.
Then $\wfl(\bZ)$ consists of the ghost vectors $\AngBrak{a_1,a_2,\dots}$ that satisfy 
	$$
	a_{j}\equiv a_{pj}\bmod p^{1+\ord_p(j)}
	$$ 
for all $j\geq 1$ and all primes $p$.

\subsection{} {\em Representing $\wfl$.}
\label{subsec:representing-wfl}
Let us construct a flat $R$-algebra $\Lambda_{R,\ptst}$
representing the functor $\wfl_{R,\ptst}$.
First we will construct objects $\Lambda^i_{R,\ptst}$ 
representing the functors $U_i$.
For $i=0$, it is clear: $U_0$ is represented by $\Lambda^0_{R,\ptst}=\Psi_{R,\ptst}$.
Now assume $\Lambda^i_{R,\ptst}$ has been constructed and that
it is a sub-$R$-algebra of $R[1/\ptst]\tn_R \Psi_{R,\ptst}$ satisfying
	$$
	R[1/\ptst]\tn_R\Lambda^i_{R,\ptst}=R[1/\ptst]\tn_R \Psi_{R,\ptst}.
	$$ 
Then let $\Lambda^{i+1}_{R,\ptst}$ denote the sub-$\Lambda^i_{R,\ptst}$-algebra
of $R[1/\ptst]\tn_R \Psi_{R,\ptst}$ generated by 
all elements $\pi^*\tn(\psi_\pind(f)-f^{q_\pind})$,
where $\pi^*\in\m_\pind^{-1}\subseteq R[1/\ptst]$,
$f\in\Lambda^i_{R,\ptst}$, and $\pind\in\ptst$.
Then $\Lambda^i_{R,\ptst}$ is flat over $R$.  Indeed, it is $\ptst$-flat
because it is a sub-$R$-algebra of $R[1/\ptst]\tn_R \Psi_{R,\ptst}$, 
and it is flat away from $\ptst$ because $R[1/\ptst]\tn_R\Lambda_{R,E}$
agrees with the free $R[1/\ptst]$-algebra $R[1/E]\tn_R\Psi_{R,\ptst}$.
It also clearly represents $U_i$. 

Finally, we set
	\begin{equation}
		\Lambda_{R,\ptst} = \bigcup_{i\geq 0} \Lambda^i_{R,\ptst} \subseteq 
			R[1/\ptst]\tn_R\Psi_{R,\ptst}.
	\end{equation}
It is flat over $R$ because it is a colimit of flat $R$-algebras,
and it represents $\wfl_{R,\ptst}$ because each $\Lambda^i_{R,\ptst}$ represents $U_i$. 
As an example, if $\ptst=\ptst'\smcoprod\ptst''$, where $\ptst''$ consists 
of only copies of the unit ideal, then $\Lambda_{R,\ptst}$ agrees with the
monoid algebra $\Lambda_{R,\ptst'}[\bN^{(\ptst'')}]$.
We will often use the shortened forms $\Lambda_{\ptst}$ or, when $\ptst=\{\m\}$,
$\Lambda_{\m}$.

Since $\Lambda_{R,\ptst}$ represents $\wfl$, which takes
values in $R$-algebras, $\Lambda_{R,\ptst}$ carries
the structure of a co-$R$-algebra object in $\Ringfl_R$.
Because $\Ringfl_R$ is closed under coproducts (the tensor product of flat modules
being flat), a co-ring structure consists in morphisms
	\begin{equation} \label{maps:biring-structure-on-Lambda}
	\Delta^+,\Delta^{\times}\:\Lambda_{R,\ptst}\longmap \Lambda_{R,\ptst}\tn_R \Lambda_{R,\ptst},
	\quad \eps^+,\eps^{\times}\:\Lambda_{R,\ptst} \longmap R
	\end{equation}
corresponding to addition, multiplication, the additive identity, and the multiplicative identity
on the functor $\wfl_{R,\ptst}$. The $R$-linear structure on $\wfl_{R,\ptst}$ corresponds to a
morphism
	\begin{equation} \label{maps:co-R-structure-on-Lambda}
		\beta\:\Lambda_{R,\ptst} \to R^R = \prod_R R.
	\end{equation}
All these structure maps satisfy the opposite of the $R$-algebra axioms.
(In the language of schemes, one would say this makes $\Spec \Lambda_{R,\ptst}$ an
$R$-algebra scheme over $R$; or in the language of~\cite{Borger-Wieland:PA}, it makes
$\Lambda_{R,\ptst}$ an $R$-$R$-biring.)

\subsection{} {\em Definition of $W$ in general.}
\label{subsec:W-define-general}
We can view $\Lambda_{R,\ptst}$ as an object of $\Ring_R$, instead of $\Ringfl_R$.
Then define $W_{R,\ptst}$ as a set-valued functor on $\Ring_R$ by
	\begin{equation}
		W_{R,\ptst}(A) = \Hom_{\Ring_R}(\Lambda_{R,\ptst},A).
	\end{equation}
The structure maps (\ref{maps:biring-structure-on-Lambda})--(\ref{maps:co-R-structure-on-Lambda})
give $W_{R,\ptst}$ the structure of a functor with values in $R$-algebras:
	\begin{equation}
		W_{R,\ptst}\:\Ring_R\longmap\Ring_R.
	\end{equation}
(Note that here we really use the fact the the coproduct in $\Ringfl_R$ agrees with that
in $\Ring_R$. In \ref{subsec:representing-wfl}, it was used only to justify the
symbol $\tn$ for the coproduct.)

For any $A\in\Ring_R$, let us call the $W_{R,\ptst}(A)$ the $R$-algebra of \emph{$\ptst$-typical
Witt vectors with entries in $A$}. Its restriction to $\Ringfl_R$ agrees with $\wfl_{R,\ptst}$
because $\Ringfl_R$ is a full subcategory of $\Ring_R$.

We will often write $W_{\ptst}$ or $W$ for $W_{R,\ptst}$ when there is no risk of confusion.
When $\ptst$ consists of a single ideal $\m$, we will also write $W_{R,\m}$ or $W_{\m}$.

\subsection{} \emph{Remark: Kan extensions.}
\label{subsec:Kan-extensions}
In categorical terms, 
$W_{R,\ptst}$ is the left Kan extension of $i\circ\wfl_{R,\ptst}$ along the inclusion
functor $i$:
	\begin{equation} \label{diag:W-Kan-ext}
		\xymatrix{
		\Ringfl_R \ar^i[r]		
			& \Ring_R \\
		\Ringfl_R \ar^{\wfl_{R,\ptst}}[u]\ar^i[r]	
			& \Ring_R. \ar_{W_{R,\ptst}}@{-->}[u]
		}
	\end{equation}
(See Borceux~\cite{Borceux:Handbook.v1}, 3.7, for example, for the general theory of Kan
extensions.) I mention this only to emphasize that the passage from the $\ptst$-flat case to the
general case is by a purely category-theoretic process, and hence the heart of the theory lies in
the $\ptst$-flat case. This is in contrast to the common point of view that the purpose of
Witt vector functors is to lift rings from positive characteristic to characteristic zero.

\subsection{} {\em Ghost map $\gh{}$.}
\label{subsec:def-ghost-map-for-rings}
The ghost map
	\[
	\gh{}\: W_{R,\ptst}(A) \longmap \prod_{\nset} A
	\]
is the natural map induced by the universal property of Kan extensions
applied to the inclusion maps $\wfl_{R,\ptst}(A) \to \prod_{\nset}A$, 
which are functorial in $A$.  Equivalently, it is the morphism of functors induced
by the map
	$$
	\Psi_{R,\ptst} = \Lambda^0_{R,\ptst} \longmap \Lambda_{R,\ptst}
	$$
of representing objects.  When $A$ is $\ptst$-flat, it is harmless to 
identify $\gh{}$ with the inclusion map.

\subsection{} {\em Example: $p$-typical and big Witt vectors.}
\label{subsec:agreement-with-classical-witt-vectors}
Suppose $R$ is $\bZ$.
If $\ptst$ consists of the single ideal $p\bZ$, then $W$ agrees with
the classical $p$-typical Witt vector functor~\cite{Witt:Vectors}.
Indeed, for $p$-torsion free rings $A$, this follows from Cartier's lemma, which says that
the traditionally defined $p$-typical Witt vector functor restricted to the category
of $p$-torsion-free rings has the same universal property as
$\wfl$. (See Bourbaki~\cite{Bourbaki:CommAlg}, IX.44, exercise 14 or Lazard~\cite{Lazard:CFG},
VII\S 4.)  Therefore, they are isomorphic functors.  For $A$ general, one just observes
that that the traditional functor is represented by the ring $\bZ[x_0,x_1,\dots]$,
which is $p$-torsion free, and so it is the left Kan extension of its
restriction to the category of $p$-torsion-free rings.  Therefore it agrees with $W$
as defined here.

Another proof of this is given in~\ref{subsec:Witt-components}.  It makes a direct connection with 
the traditional Witt components, rather than going through the universal property.

Suppose instead that $\ptst$ is the family of all maximal ideals of $\bZ$.
Then $W$ agrees with the classical big Witt vector functor.
As above, this can be shown by reducing to the torsion-free case and then
citing the analogue of Cartier's lemma.  (Which version of Cartier's lemma
depends on how we define the classical big Witt vector functor.
If we use generalized Witt polynomials,
then we need Bourbaki~\cite{Bourbaki:CommAlg}, IX.55, exercise 41b.
If it is defined as the cofree $\lambda$-ring functor, as in
Grothendieck~\cite{Grothendieck:Chern}, then
we need Wilkerson's theorem~\cite{Wilkerson}, proposition 1.2.)

Finally, we will see below (\ref{subsec:Witt-components}) that when $R$ is a complete discrete
valuation ring and $\ptst$ consists of the maximal ideal of $R$, then $W$ agrees with
Hazewinkel's ramified Witt vector functor~\cite{Hazewinkel:book}, (18.6.13).

\subsection{} {\em Comonad structure on $W$.}
The functor $\wfl\:\Ringfl_R\to\Ringfl_R$ is naturally a comonad, being the composition of a
functor (the forgetful one) with its right adjoint, and this comonad structure prolongs naturally
to $W_{R,\ptst}$. The reason for this can be expressed in two ways---in terms of Kan extensions or
in terms of representing objects.

The first way is to invoke the general fact that $W_{R,\ptst}$, as the Kan extension of the comonad
$\wfl_{R,\ptst}$, has a natural comonad structure. This uses the commutativity
of~(\ref{diag:W-Kan-ext}) and the fullness and faithfulness of $i$. The other way is to translate
the structure on $\wfl$ of being a comonad into a structure on its representing object
$\Lambda_{R,\ptst}$. One then observes that this is exactly the structure for the underlying
$R$-algebra $i(\Lambda_{R,\ptst})$ to represent a comonad on $\Ring_R$. (This is called an
$R$-plethory structure in~\cite{Borger-Wieland:PA}.)

\subsection{} {\em $\Lambda$-rings.}
\label{subsec:def-lambda-ring}
The category $\Ring_{\Lambda_{R,\ptst}}$ of $\Lambda_{R,\ptst}$-rings
is by definition the category of coalgebras for the comonad $W_{R,\ptst}$, that is,
the category of $R$-algebras equipped with a coaction of the comonad
$W_{R,\ptst}$.  Since $W_{R,\ptst}$ extends $\wfl_{R,\ptst}$,
a $\Lambda_{R,\ptst}$-ring structure on an $\ptst$-flat $R$-algebra
$A$ is the same as a commuting family of Frobenius lifts $\psi_\pind$.

When $R=\bZ$ and $\ptst$ is the family of all maximal ideals of $\bZ$, then a $\Lambda$-ring is the
same as a $\lambda$-ring in the sense of Grothendieck's Riemann--Roch
theory~\cite{Grothendieck:Chern} (and originally called a ``special $\lambda$-ring''). In the
$\ptst$-flat case, this is Wilkerson's theorem~(\cite{Wilkerson}, proposition 1.2). The proof is an
exercise in symmetric functions, but the deeper meaning eludes me.  The general case
follows from the $\ptst$-flat case by category theory, as
in~\ref{subsec:agreement-with-classical-witt-vectors}.

\subsection{} {\em Free $\Lambda$-rings and $\Lambda\bcp\vbl$.}
\label{subsec:free-lambda-rings}
Since $W_{R,\ptst}$ is a representable comonad on $\Ring_R$, the forgetful
functor from the category of $\Lambda_{R,\ptst}$-rings to the category
of $R$-algebras has a left adjoint denoted $\Lambda_{R,\ptst}\bcp\vbl$.
This follows either from the adjoint functor theorem in category 
theory (3.3.3 of Borceux~\cite{Borceux:Handbook.v1}),
or by simply writing down the adjoint in terms of
generators and relations, as in 1.3 of Borger--Wieland~\cite{Borger-Wieland:PA}.
The second approach involves the $R$-plethory structure on $\Lambda_{R,\ptst}$, and is
similar to the description of tensor products, free differential rings, and so on
in terms of generators and relations.

The functor $\Lambda_{R,\ptst}\bcp\vbl$, viewed as an endofunctor on the category
of $R$-algebras, is naturally a monad, simply because it is the left adjoint of the comonad
$W_{R,\ptst}$.  Further, the category of algebras for this monad is naturally
equivalent to $\Ring_{\Lambda_{R,\ptst}}$.
This can be proved using Beck's theorem (4.4.4 of Borceux~\cite{Borceux:Handbook.v2}), and is the 
same as the fact that the category of $K$-modules, for any ring 
$K$, can be defined as the category of algebras for the monad $K\tn\vbl$ or coalgebras for
the comonad $\Hom(K,\vbl)$.

We can interpret elements of $\Lambda_{R,\ptst}$ as natural operations on 
$\Lambda_{R,\ptst}$-rings.
Indeed, a $\Lambda_{R,\ptst}$-ring structure on a ring $A$ is by definition a (type of)
map $A\to W_{R,\ptst}(A)$. It therefore induces a set map
	$$
	\Lambda_{R,\ptst} \times A \longmap \Lambda_{R,\ptst}\times W_{R,\ptst}(A)
	 	= \Lambda_{R,\ptst} \times \Hom_R(\Lambda_{R,\ptst},A) \longmap A,
	$$
which is functorial in $A$. In particular, if we take $A=\Lambda_{R,\ptst}$, we get a 
set map
	\begin{equation} \label{map:composition-on-Lambda}
		\Lambda_{R,\ptst}\times \Lambda_{R,\ptst} \longlabelmap{\circ} \Lambda_{R,\ptst}.
	\end{equation}
It agrees with the restriction of the composition map $\circ$ on 
$\Psi_{R[1/\ptst],\ptst}=R[1/\ptst]\tn_R \Psi_{R,\ptst}$
given in~\ref{subsec:Psi-rings}. In particular, it is associative with identity $e$.

In fact, all natural operations on $\Lambda_{R,\ptst}$-rings come from
$\Lambda_{R,\ptst}$ in this way.  See \cite{Borger-Wieland:PA} for an abstract
account from this point of view.

\subsection{} {\em Remark: identity-based approaches.}
\label{rmk:formula-based-approaches}
It is possible to set up the theory of $\Lambda_{R,\ptst}$-rings
more concretely using
universal identities rather than category theory. (See Buium, Buium--Simanca, and
Joyal
\cite{Buium:p-jets}\cite{Buium-Simanca:Arithmetic-Laplacians}\cite{Joyal:Lambda}\cite{Joyal:Witt},
for example.) In this subsection,
I will say something about that point of view and its relation to the
category-theoretic one, but it will not be used elsewhere in this paper.

First suppose that for each $\pind\in\ptst$, the ideal $\m_\pind$ is generated by a 
single element $\pi_\pind$.  For any $\Lambda_{R,\ptst}$-ring $A$
and any element $a\in A$, there exists an element $\delta_{\pind}(a)\in A$ such that
	$$
	\psi_{\pind}(a) = a^{q_{\pind}} + \pi_\pind \delta_{\pind}(a).
	$$
If we now assume that $A$ is $\ptst$-flat, then the element $\delta_{\pind}(a)$ is uniquely
determined by this equation,
and therefore $\delta_{\pind}$ defines an operator on $A$:
	$$
	\delta_{\pind}(a) = \frac{\psi_{\pind}(a)-a^{q_{\pind}}}{\pi_{\pind}}.
	$$
Observe that if the integer $q_{\alpha}$ maps to $0$ in $R$, for example when $R$
is a ring of integers in a function field,
then $\delta_{\pind}$ is additive; but otherwise it essentially never is.
(Also note that $\delta_{\pind}$ is the same as the operator $\theta_{\pi_\pind,1}$ defined 
in~\ref{subsec:theta-operator-def} below.)

Conversely, if $A$ is an $\ptst$-flat $R$-algebra, equipped with with operators
$\delta_\pind$, then there is at most one $\Lambda_{R,\ptst}$-ring structure on $A$
whose $\delta_{\pind}$-operators are the given ones.  To say when such a 
$\Lambda_{R,\ptst}$-ring structure exists, we only need to express
in terms of the operators $\delta_{\pind}$ the condition
that the operators $\psi_{\pind}$ be commuting $R$-algebra homomorphisms.  After 
dividing by any accumulated factors of $\pi_{\pind}$, this gives the 
identities of Buium--Simanca~\cite{Buium-Simanca:Arithmetic-Laplacians},
definition 2.1:  
	\begin{align}
		\label{eq:delta-ring-axiom-1}
		\delta_{\pind}(r) &= \frac{r-r^{q_{\pind}}}{\pi_{\pind}}, \quad\text{for } r\in R,\\
		\label{eq:delta-ring-axiom-2}
		\delta_{\pind}(a+b) &= \delta_{\pind}(a) + \delta_{\pind}(b) + C_{\pind}(a,b), \\
		\label{eq:delta-ring-axiom-3}
		\delta_{\pind}(ab) &= \delta_{\pind}(a)b^{q_{\pind}} + a^{q_{\pind}}\delta_{\pind}(b)
			+ \pi_{\pind} \delta_{\pind}(a)\delta_{\pind}(b), \\
		\label{eq:delta-ring-axiom-4}
		\delta_{\pind}\circ\delta_{\pind'}(a) &= \delta_{\pind'}\circ\delta_{\pind}(a)
			+ C_{\pind,\pind'}\big(a,\delta_{\pind}(a),\delta_{\pind'}(a)\big),
	\end{align}
where 
\begin{equation}
		\label{eq:delta-ring-axiom-5}
		C_{\pind}(x,y) = \frac{x^{q_\pind} + y^{q_{\pind}}-(x+y)^{q_{\pind}}}{\pi_{\pind}}
			=-\sum_{i=1}^{q_\pind-1}\frac{1}{\pi_\pind}\binom{q_\pind}{i}x^{q_\pind-i}y^i, \\
\end{equation}
	\begin{multline}
		\label{eq:delta-ring-axiom-6}
		C_{\pind,\pind'}(x,y,z) = 
			\frac{C_{\pind'}(x^{q_{\pind}},\pi_{\pind}y)}{\pi_{\pind}}
			- \frac{C_{\pind}(x^{q_{\pind'}},\pi_{\pind'}z)}{\pi_{\pind'}} \\
			- \frac{\delta_{\pind}(\pi_{\pind'})}{\pi_{\pind'}}z^{q_{\pind}} 
			+ \frac{\delta_{\pind'}(\pi_{\pind})}{\pi_{\pind}}y^{q_{\pind'}}.
	\end{multline}
One can easily check that the coefficients of these polynomials are elements of
$R$.

For any $R$-algebra $A$, let us define a $\delta_{R,\ptst}$-structure on $A$ to be a family of
operators $\delta_{\pind}$ satisfying the axioms above. Thus, if $A$ is an $\ptst$-flat
$R$-algebra, then a $\Lambda_{R,\ptst}$-structure---by definition a commuting family of Frobenius
lifts indexed by $\ptst$---is equivalent to a $\delta_{R,\ptst}$-structure. The point of all this,
then, is that if we no longer require $A$ to be $\ptst$-flat, a $\delta_{R,\ptst}$-structure is
generally stronger than having a commuting family of Frobenius lifts, but it is still equivalent to
having a $\Lambda_{R,\ptst}$-structure. This offers another point of view on the difference
between a $\Lambda_{R,\ptst}$-structure and a commuting family of Frobenius lifts: A
$\delta_{R,\ptst}$-structure is well behaved from the point of view of universal algebra (and hence
so is a $\Lambda_{R,\ptst}$-structure) because it is given by operators $\delta_{\pind}$ whose
effect on the ring structure is described by universal identities, as above; but the structure of a
commuting family of Frobenius lifts does not have this property because of the existential
quantifier hidden in the word \emph{lift}.

The equivalence between $\delta_{R,\ptst}$-structures and $\Lambda_{R,\ptst}$-structures can be
seen as follows. For $\ptst$-flat $R$-algebras $A$, it was explained above. For general $A$, the
equivalence can be shown by checking that the cofree $\delta_{R,\ptst}$-ring functor is represented
by an $\ptst$-flat $R$-algebra (in fact, a free one). It therefore agrees with the left Kan
extension of its restriction to the category of $\ptst$-flat algebras, and hence agrees with
$W_{R,\ptst}$.

We could extend the identity-based approach to the case where
the ideals $\m_\pind$ are not principal, but then we would need operators 
	\begin{equation} \label{eq:theta-def}
		\delta_{\pind,\pi_{\pind}^*}(x) = \pi_{\pind}^*(\psi_\pind(x)-x^{q_\pind})
	\end{equation}
for every element $\pi_{\pind}^*\in\m_\pind^{-1}$, or at least for 
those in a chosen generating set of $\m_{\pind}^{-1}$, and
we would need additional axioms relating them. A particularly convenient generating set of 
$\m_{\pind}^{-1}$
is one of the form $\{1,\pi_{\pind}^*\}$, which always exists. Further, for each $\pind\in\ptst$,
it is enough to use the operators $\psi_{\pind}$ and $\delta_{\pind,\pi_{\pind}^*}$
instead of $\delta_{\pind,1}$ and $\delta_{\pind,\pi_{\pind}^*}$,
because $\delta_{\pind,1}$ can be expressed in terms of $\psi_{\pind}$, by~(\ref{eq:theta-def}).
Therefore if we fix elements $\pi_{\pind}^*\in \m_{\pind}^{-1}$ which are $R$-module
generators modulo $1$,
the relations needed for the generating set
$\bigcup_{\pind\in\ptst}\{\psi_{\pind},\delta_{\pind,\pi_{\pind}^*}\}$ of operators
are those in (\ref{eq:delta-ring-axiom-1})--(\ref{eq:delta-ring-axiom-6}) 
but one needs to make the following changes: for each $\pind\in\ptst$,
replace each occurrence of $\pi_{\pind}^{-1}$ with $\pi_{\pind}^*$, and add axioms that
$\psi_{\pind}$ is
an $R$-algebra homomorphism, 
that $\psi_{\pind}$ commutes with all $\psi_{\pind'}$ and all $\delta_{\pind',\pi_{\pind'}^*}$,
and that (\ref{eq:theta-def}) holds.

When $R$ is an $\bF_p$-algebra for some prime number $p$, the polynomials
$C_{\pind}(x,y)$ are zero and the axioms above simplify considerably. In particular, the
operators $\delta_{\pind}$ are additive, and so it is possible to describe a
$\Lambda_{R,\ptst}$-structure using a cocommutative twisted bialgebra, the additive bialgebra of
the plethory $\Lambda_{R,\ptst}$. (See Borger--Wieland~\cite{Borger-Wieland:PA}, sections 2 and
10.)

\subsection{} \emph{Localization of the ring $R$ of scalars.}
\label{subsec:localizing-R-infinite-length}
Let $R'$ be an $\ptst$-flat $R$-algebra such that the structure map $R\to R'$ is an epimorphism of
rings. (For example, the map $\Spec R'\to\Spec R$ could be an open immersion.) Then the family
$(\m_{\pind})_{\pind\in\ptst}$ induces a family $(\m'_{\pind})_{\pind\in\ptst}$ of ideals of $R'$,
where $\m'_{\pind}=\m_{\pind}R'$. By the assumptions on $R'$, each $\m'_{\pind}$ is supramaximal.
Let us write $\ptst'=\ptst$ and use the notation $\ptst'$ for the index set of the $\m'_{\pind}$.

Let us construct an isomorphism:
	\begin{equation} \label{map:lambda-change-of-R-infinite-length}
		R'\tn_R\Lambda_{R,\ptst} \longisomap \Lambda_{R',\ptst'}.
	\end{equation}

The category $\Ringfl_{\Lambda_{R',\ptst'}}$ (see~\ref{subsec:E-flat-lambda-rings}) is a
subcategory of the category of $\Ringfl_{\Lambda_{R,\ptst}}$. Indeed, any object
$A'\in\Ringfl_{\Lambda_{R',\ptst'}}$ is an $R$-algebra with endomorphisms $\psi_{\m_{\pind}}$, for
each $\pind\in\ptst$. These endomorphisms are again commuting Frobenius lifts, simply because
$A'/\m'_{\pind}A'=A'/\m_{\pind}A'$. Since $A'$ is $\ptst'$-flat (and by the assumptions on $R'$),
$A'$ is $\ptst$-flat. Therefore, it can be viewed as a $\Lambda_{R,\ptst}$-ring.

Further, $\Ringfl_{\Lambda_{R',\ptst'}}$ agrees with the subcategory of
$\Ringfl_{\Lambda_{R,\ptst}}$ consisting of objects $A$ whose structure map $R\to A$ factors
through $R'$, necessarily uniquely. Now consider the underlying-set functor on this category. From
the definition of $\Ringfl_{\Lambda_{R',\ptst'}}$, this functor is represented by the right-hand
side of~(\ref{map:lambda-change-of-R-infinite-length}), and from the second description, it is
represented by the left-hand side. Let~(\ref{map:lambda-change-of-R-infinite-length}) be the
induced isomorphism on representing objects. It sends an element $r'\tn f$ to $r'f.$

The isomorphism of represented functors which is
induced by~(\ref{map:lambda-change-of-R-infinite-length}) 
gives natural maps 
	\begin{equation}\label{map:witt-independence-of-R-infinite-length}
		W_{R',\ptst'}(A') \longisomap W_{R,\ptst}(A'),
	\end{equation}
for $R'$-algebras $A'$.

Finally, let us show that for any $R'$-algebra $B'$,
the following canonical map is an isomorphism:
	\begin{equation} \label{eq:lambda-circle-independent-of-R-infinite-length}
		\Lambda_{R,\ptst}\bcp B' \longisomap \Lambda_{R',\ptst'}\bcp B'.
	\end{equation}
It is enough to show that for any $R'$-algebra $A'$, the induced map
	$$
	\Hom_{R'}(\Lambda_{R',\ptst'}\bcp B',A')\longmap\Hom_{R'}(\Lambda_{R,\ptst}\bcp B', A')
	$$
is a bijection.
Since $\Ring_{R'}$ is a full subcategory of $\Ring_R$, 
the right-hand side agrees with $\Hom_{R}(\Lambda_{R,\ptst}\bcp B', A')$,
and so the map above is an isomorphism by~(\ref{map:witt-independence-of-R-infinite-length}).

\subsection{} {\em Teichm\"uller lifts.}
\label{subsec:teich-1}
Let $A$ be an $R$-algebra, let $A^\circ$ denote the commutative monoid of all elements of $A$ 
under multiplication, and let $R[A^{\circ}]$ denote the monoid algebra on 
$A^\circ$.  Then for each $\pind\in\ptst$,
the monoid endomorphism $a\mapsto a^{q_\pind}$ of $A^\circ$ induces an $R$-algebra 
endomorphism $\psi_\pind$ of $R[A^{\circ}]$ which reduces to the $q_\pind$-th power
map modulo $\m_\pind$.  Since $R[A^\circ]$ is free as an $R$-module,
it is flat.  And since the various $\psi_\pind$ commute with each other,
they provide $R[A^\circ]$ with a $\Lambda_{R,\ptst}$-structure.
Combined with the $R$-algebra map $R[A^\circ]\to A$ given by the counit of the evident 
adjunction, this gives, by the right-adjoint property of $W_{R,\ptst}$,
a $\Lambda_{R,\ptst}$-ring map $t\:R[A^\circ]\to W_{R,\ptst}(A)$.  We write the composite 
monoid map
	\[
	A^\circ \longlabelmap{\mathrm{unit}} R[A^\circ]^\circ \longlabelmap{t^\circ} 
		W_{R,\ptst}(A)^\circ
	\]
as simply $a\mapsto [a]$.  It is a section of the $R$-algebra map 
$\gh{0}\:W_{R,\ptst}(A)\to A$
and is easily seen to be functorial in $A$.
The element $[a]$ is called the {\em Teichm\"uller lift} of $a$.

\section{Grading and truncations}
\label{sec:gradings-and-truncations}

\subsection{} \emph{Ordering on $\zset$.}
For two elements $n',n\in\zset=\bigoplus_{\ptst}\bZ$,
write $n'\leq n$ if we have $n'_\pind\leq n_\pind$ for all $\pind\in\ptst$.
Also put 
	$$
	[0,n] = \setof{n'\in\nset}{n'\leq n}.
	$$

\subsection{} \label{subsec:grading-and-truncations}
{\em Truncations.}
We have the following decomposition of $\Psi_{R,\ptst}$:
\[
\Psi_{R,\ptst}=\bigotimes_{\pind\in \ptst}\bigotimes_{i\in\bN}R[\psi_\pind^{\circ i}]\,
	= \bigotimes_{n\in\nset} R[\psi_n] = R[\psi_n | n\in\nset].
\]
(Thus, $\Psi_{R,\ptst}$ is an $\nset$-indexed coproduct in the category of 
$R$-algebras,
much like graded rings are monoid-indexed coproducts in the category of modules. One might say
that $\Psi_{R,\ptst}$ is an $\nset$-graded plethory. This point of view will not 
be used below.)
For each $n\in\zset$, put
	\[
	\Psi_{R,\ptst,n} = \bigotimes_{\pind\in \ptst}\,\bigotimes_{0\leq i\leq 
	n_\pind}R[\psi_\pind^{\circ i}] = \bigotimes_{n'\in [0,n]}R[\psi_{n'}] 
	= R\big[\psi_{n'}|n'\in [0,n]\big].
	\]
Then $\Psi_{R,\ptst,n}$ represents the $\Ring_R$-valued functor that sends $A$ to
the product ring $A^{[0,n]}$, which is naturally a quotient of $A^{\nset}$. 

Define a similar filtration on $\Lambda_{R,\ptst}$ by
	\begin{equation} \label{eq:filtration-on-Lambda}
		\Lambda_{R,\ptst,n} = \Lambda_{R,\ptst}\cap\big(R[1/\ptst]\tn_R\Psi_{R,\ptst,n}\big).
	\end{equation}
We will often use the shortened forms $\Lambda_{\ptst,n}$, $\Psi_{\ptst,n}$,
$\Lambda_{\m,n}$, $\Psi_{\m,n}$, and so on.

\begin{proposition}\label{pro:degree-filt-on-Lambda}
	\begin{enumerate}
		\item For each $n\in\nset$, the $R$-scheme $\Spec\Lambda_{R,\ptst,n}$ admits
			a unique structure of an $R$-algebra object in the category of $R$-schemes
			such that the map $\Spec\Lambda_{R,\ptst}\to \Spec\Lambda_{R,\ptst,n}$ 
			induced by the inclusion $\Lambda_{R,\ptst,n}\subseteq\Lambda_{R,\ptst}$
			is a homomorphism of $R$-algebra schemes over $R$.
		\item For each $m,n\in\nset$, we have
			\begin{equation} \label{eq:plethysm-def}
			\Lambda_{R,\ptst,m}\circ\Lambda_{R,\ptst,n}\subseteq\Lambda_{R,\ptst,m+n},
			\end{equation}
			where $\circ$ denotes the composition map
			of (\ref{map:composition-on-Lambda}).
	\end{enumerate}
\end{proposition}

\begin{proof}
	(a): 
	Write $\Lambda=\Lambda_{R,\ptst}$, $\Lambda_n=\Lambda_{R,\ptst,n}$, and so on.
	First observe that, for any integer $i\geq 0$, all the maps in the diagram
		$$
		\xymatrix{
		R[1/\ptst]\tn_R \Psi_n^{\tn_R i} \ar^{a_i}[r] & R[1/\ptst]\tn_R \Psi^{\tn_R i} \\
		\Lambda_n^{\tn_R i} \ar^{b_i}[r]\ar[u] & \Lambda^{\tn_R i} \ar^{c_i}[u] 
		}
		$$
	are injective. Indeed, $a_i$ clearly is; the vertical maps are because they become isomorphisms
	after base change to $R[1/\ptst]$ and because $\Lambda_n$ and $\Lambda$
	are $\ptst$-flat; and it follows formally that $b_i$ is injective.
	Then the uniqueness of the desired $R$-algebra scheme structure on $\Spec \Lambda_n$,
	follows from the injectivity of $b_2$.

	Now consider existence. Let
		$$
		\Delta\:R[1/\ptst]\tn_R \Psi \longmap R[1/\ptst]\tn_R \Psi\tn_R\Psi
		$$
	denote the ring map that induces the addition (resp.\ multiplication) map
	on the ring scheme $\Spec R[1/\ptst]\tn_R \Psi_R$.
	To show that the desired addition and multiplication maps on $\Spec\Lambda_n$ exist,
	it is enough to show 
		\begin{equation} \label{eq:internal-hiho}
			\Delta(\Lambda_n)\subseteq \Lambda_n\tn_R\Lambda_n.
		\end{equation}
	In fact, once we do this, we will be done: because each $c_i\circ b_i$ is injective, the ring
	axioms (associativity, distributivity, \dots) will follow from those on
	$\Spec R[1/\ptst]\tn_R\Psi$.

	The map $\Delta$ sends $\psi_\pind$ to
	$\psi_{\pind}\tn 1 + 1\tn \psi_{\pind}$ (resp. $\psi_{\pind}\tn \psi_{\pind}$).	
	Therefore we have
		$$
		\Delta(R[1/\ptst]\tn_R \Psi_n) \subseteq R[1/\ptst]\tn_R \Psi_n\tn_R\Psi_n,
		$$
	and hence
		$$
		\Delta(\Lambda_{n}) \subseteq \Lambda^{\tn_R 2} 
			\cap \big(R[1/\ptst]\tn_R \Psi_{n}^{\tn_R 2}\big) = \Lambda_{n}^{\tn_R 2}.
		$$
	This establishes (\ref{eq:internal-hiho}) and hence completes the proof of (a).

	(b): Combine the definition~(\ref{eq:filtration-on-Lambda}) with the inclusion
		\[
		(R[1/\ptst]\tn_R\Psi_{m})\circ(R[1/\ptst]\tn_R\Psi_{n})
			\subseteq(R[1/\ptst]\tn_R\Psi_{m+n})
		\]
	and the inclusion
	$\Lambda_{m}\circ\Lambda_{n}\subseteq \Lambda$.
\end{proof}

\subsection{} {\em Witt vectors of finite length.}
\label{subsec:W-finite-length}
Let $W_{R,\ptst,n}$ denote the functor $\Ring_R\to\Ring_R$
represented by $\Lambda_{R,\ptst,n}$:
	\begin{equation} \label{eq:finite-length-witt-vector-def}
		W_{R,\ptst,n}(A)=\Hom_R(\Lambda_{R,\ptst,n},A).
	\end{equation}
We call $W_{R,\ptst,n}$ the \emph{$\ptst$-typical Witt vector functor of length $n$}.
As in \ref{subsec:W-define-general},  we will often write $W_{\ptst,n}$ or $W_n$; when $\ptst=\{\m\}$, we
will also write $W_{R,\m,n}$ or $W_{\m,n}$.
We then have
	\begin{equation} \label{eq:W-inf-is-inverse-limit-of-Wn}
		W_{R,\ptst}(A) = \limm_n W_{R,\ptst,n}(A).
	\end{equation}
(Note that it is often
better to view $W_{R,\ptst}(A)$ as a pro-ring than to actually take the limit.
If we preferred topological rings to pro-rings, we could 
take the limit and endow it with the natural pro-discrete topology.) It follows
from \ref{pro:V-exact-sequence} and (\ref{map:witt-iterate-generalized-finite-length}) below that 
the maps in this projective system are surjective.

The (truncated) ghost map
	\begin{equation} \label{map:ghost-map-for-rings-leq-n}
		\gh{\leq n} \: W_{R,\ptst,n}(A) \longmap  A^{[0,n]},	
	\end{equation}
is the one induced by the inclusion $\Psi_{R,\ptst,n}\subseteq \Lambda_{R,\ptst,n}$
of representing objects.
For any $i\in[0,n]$, the composition $\gh{\leq n}$ with the 
projection onto the $i$-th factor gives another natural map
	\begin{equation} \label{map:ghost-map-for-rings-equal-n}
		\gh{i}  \: W_{R,\ptst,n}(A) \longmap A.
	\end{equation}
Also the containment~(\ref{eq:plethysm-def}) induces an $R$-algebra map
	\begin{equation} \label{map:witt-coplethysm-finite-length}
		W_{R,\ptst,m+n}(A)\longmap W_{R,\ptst,n}\big(W_{R,\ptst,m}(A)\big)
	\end{equation}
which sends an element $a\:\Lambda_{R,\ptst,m+n}\to A$ of $W_{R,\ptst,m+n}(A)$ to the map
$\gamma\mapsto [\beta\mapsto a(\beta\circ\gamma)]$, for variables $\gamma\in\Lambda_{R,\ptst,n}$
and $\beta\in\Lambda_{R,\ptst,m}$. We will call (\ref{map:witt-coplethysm-finite-length}) 
\emph{co-plethysm}.
It agrees with the map of functors induced by the map
	\begin{equation}
		\Lambda_{R,\ptst,m}\bcp\Lambda_{R,\ptst,n}\longmap \Lambda_{R,\ptst,m+n}, \quad
		\beta\bcp\gamma\mapsto\beta\circ\gamma
	\end{equation}
on representing objects, where $\beta\bcp\gamma$ is defined as in \cite{Borger-Wieland:PA}.

Finally, observe that for any element $f\in\Lambda_{R,\ptst,n}$ the natural 
$\Lambda_{R,\ptst}$-ring operation 
	$$
	f\:W_{R,\ptst}(A)\to W_{R,\ptst}(A)
	$$
(a map of sets) descends to a map $f\:W_{R,\ptst,m+n}(A)\to W_{R,\ptst,m}(A)$.
Indeed, it is the composition
	\begin{multline}
		W_{R,\ptst,m+n}(A)\longlabelmap{(\ref{map:witt-coplethysm-finite-length})} 
			W_{R,\ptst,n}\big(W_{R,\ptst,m}(A)\big)
			= \Hom(\Lambda_{R,\ptst,n},W_{R,\ptst,m}(A)) \\
			\longlabelmap{\vbl(f)} W_{R,\ptst,m}(A),
	\end{multline}
where $\vbl(f)$ denotes the map that evaluates at $f$. Particularly important
is the example $f=\psi_n$, where the induced map
	\begin{equation}
		\label{map:psi-for-finite-length-W}
		\psi_n\:W_{R,\ptst,m+n}(A)\to W_{R,\ptst,m}(A)
	\end{equation}
is a ring homomorphism.

\subsection{} \emph{Remark: traditional versus normalized indexing.}
\label{subsec:traditional-versus-normalized-indexing}
Consider the $p$-typical Witt vectors, where $R$ is $\bZ$ and $\ptst$ consists of the single ideal
$p\bZ$. Let $W'_n$ denote Witt's functor, as defined in~\cite{Witt:Vectors}. So, for example,
$W'_n(\bF_p)=\bZ/p^n\bZ$. In~\ref{subsec:Witt-components}, 
we will construct an isomorphism $W'_{n+1}\cong W_{n}$. Thus, up
to a normalization of indices, our truncated Witt functors agree with Witt's.

The reason for this normalization is to make the indexing behave well under plethysm.
By~(\ref{eq:plethysm-def}) and~(\ref{map:witt-coplethysm-finite-length}), the index set has the
structure of a commutative monoid, and so it is preferable to use an index set with a familiar
monoid structure. If we were to insist on agreement with Witt's indexing, we would have to 
replace the sum $m+n$
in~(\ref{eq:plethysm-def}) and~(\ref{map:witt-coplethysm-finite-length}) with $m+n-(1,1,\dots)$,
where this would be computed in the product group $\bZ^{\ptst}$. The reason why this has not come
up in earlier work is that the plethysm structure has traditionally been used only through the
Frobenius maps $\psi_\pind$. In other words, only the shift operator on the indexing set was used.
Thus the distinction between $\bN$ and $\bZ_{\geq 1}$ was not so important because the shift
operator $n\mapsto n+1$ is written the same way on both. But making the identification of 
$\bN$ and $\bZ_{\geq 1}$
a monoid isomorphism would involve the unwelcome addition law $m+n-1$ on $\bZ_{\geq 1}$.

It is different with the big Witt vectors, where $R$ is $\bZ$ and $\ptst$ consists of all
maximal ideals (\ref{subsec:agreement-with-classical-witt-vectors}). They are also traditionally indexed by 
the positive integers (\cite{Hazewinkel:book}, (17.4.4)), but here the positive integers are
used multiplicatively rather than additively. In particular, the monoid structure that is required
is the obvious one; so the traditional indexing is in agreement with the normalized one: the big
Witt ring $W_{p^n}(A)$ (using traditional 
multiplicative indexing) is naturally isomorphic to our $p$-typical
ring $W_n(A)$ and to Witt's $W'_{n+1}(A)$.

\subsection{} \emph{Localization of the ring $R$ of scalars.}
\label{subsec:localizing-R-finite-length}
Let $R'$ be an $\ptst$-flat $R$-algebra such that the structure map $R\to R'$ is an epimorphism
of rings, as in~\ref{subsec:localizing-R-infinite-length}.

Then for each $n\in\nset$, we have
	\begin{align*}
		R'\tn_R \Lambda_{R,\ptst,n} 
			&= R'\tn_R\big(\Lambda_{R,\ptst}\cap(R[1/\ptst]\tn_R\Psi_{R,\ptst,n})\big) \\
			& \isomap \big(R'\tn_R\Lambda_{R,\ptst}\big)\cap(R'[1/\ptst]\tn_{R'}\Psi_{R',\ptst',n}).
	\end{align*}
(We only need to check that the displayed map is an isomorphism along $\ptst$, in which case it is
is true because $R'$ is $\ptst$-flat over $R$.)
By~(\ref{map:lambda-change-of-R-infinite-length}), this gives an isomorphism 
of $R'$-algebras
	\begin{equation} \label{map:lambda-change-of-R-finite-length}
		R'\tn_R\Lambda_{R,\ptst,n} \longisomap \Lambda_{R',\ptst',n}.
	\end{equation}
The induced isomorphism of represented functors gives natural maps 
	\begin{equation}\label{map:witt-independence-of-R-finite-length}
		W_{R',\ptst',n}(A') \longisomap W_{R,\ptst,n}(A'),
	\end{equation}
for $R'$-algebras $A'$.
If $A$ is an $R$-algebra, the inverse of this map induces a map
	\begin{equation} \label{map:witt-change-of-R-finite-length}
		R'\tn_R W_{R,\ptst,n}(A) \longmap W_{R',\ptst',n}(R'\tn_R A)
	\end{equation}
We will see in~\ref{pro:Zariski-localization-on-R} that this is an isomorphism.

As with~(\ref{eq:lambda-circle-independent-of-R-infinite-length}), the
map~(\ref{map:witt-independence-of-R-finite-length}) induces an isomorphism
	\begin{equation} \label{eq:lambda-circle-independent-of-R-finite-length}
		\Lambda_{R,\ptst,n}\bcp B' \longisomap \Lambda_{R',\ptst',n}\bcp B',
	\end{equation}
for any $R'$-algebra $B'$,

\begin{proposition}\label{pro:ghost-map-inj}
	Let $A$ be an $\ptst$-flat $R$-algebra.  Then the ghost map
		\[
		\gh{\leq n}\:W_{R,\ptst,n}(A) \longmap A^{[0,n]}
		\]
	is injective.  If $A$ is an $R[1/\ptst]$-algebra, it is an
	isomorphism.
\end{proposition}

Recall that the analogous facts for infinite-length Witt vectors are also true, 
either by construction (\ref{subsec:W-for-flat-rings}) or by the universal property
(\ref{pro:wfl-univ-property}).

\begin{proof}
	If every ideal in $\ptst$ is the unit ideal, then $\Lambda_{R,\ptst}=\Psi_{R,\ptst}$, and hence 
	we have $\Lambda_{R,\ptst,n}=\Psi_{R,\ptst,n}$.
	The statement about $R[1/\ptst]$-algebras then follows 
	from~(\ref{map:lambda-change-of-R-finite-length}).
	The statement about $\ptst$-flat $R$-algebras follows by considering
	the injection $A\to R[1/\ptst]\tn_R A$ and applying the previous
	case to $R[1/\ptst]\tn_R A$.
\end{proof}

\section{Principal single-prime case}
\label{sec:principal-single-prime-case}

For this section, we will restrict to the case where $\ptst$ consists of one ideal $\m$ generated
by an element $\pi$. Our purpose is to extend the classical components of Witt vectors from the
$p$-typical context (where $R$ is $\bZ$ and $\ptst$ consists of the single ideal $p\bZ$) to this
slightly more general one. The reason for this is that the Witt components are well-suited to
calculation. In the following sections, we will see how to use them, together
with~\ref{pro:principal-principle},~\ref{cor:witt-iterate-generalized},
and~\ref{pro:Zariski-localization-on-R}, to draw conclusions when $\ptst$ is general.

In fact, the usual arguments and definitions in the classical theory of Witt vectors carry over
as long as one modifies the usual Witt polynomials by replacing every $p$ in an exponent with
$q_\m$, and every $p$ in a coefficient with $\pi$. Some things, such as the Verschiebung
operator, depend on the choice of $\pi$, and others do not, such as the Verschiebung filtration.

Let $n$ denote an element of $\bN$. Let us abbreviate $\Lambda_{\m}=\Lambda_{R,\ptst}$,
$\Lambda_{\m,n}=\Lambda_{R,\ptst,n}$, $W_{\m}=W_{R,\ptst}$,
$q=q_{\m}$, $\psi=\psi_{\m}$, and so on.

\subsection{} {\em $\theta$ operators.}
\label{subsec:theta-operator-def}
Define elements $\theta_{\pi,0},\theta_{\pi,1},\dots$ of 
	$$
	R[1/\pi]\tn_R\Lambda_{\m}=R[1/\pi]\tn_R\Psi_{\m}
	$$ 
recursively by the generalized Witt polynomials
\begin{equation}
	\label{eq:witt-comp}
		\psi^{\circ n} 
			= \theta_{\pi,0}^{q^n}+\pi\theta_{\pi,1}^{q^{n-1}}+\cdots + \pi^n\theta_{\pi,n}.
\end{equation}
(Note that the exponent on the left side means iterated composition, while the
exponents on the right mean usual exponentiation, iterated multiplication.)
As in~\ref{subsec:Psi-rings}, we can view the elements $\theta_{\pi,i}$
as natural operators on $\Psi_{R[1/\pi],\m}$-rings.  
We will often write $\theta_i=\theta_{\pi,i}$ when $\pi$ is clear.

\begin{lemma}\label{lem:psi-of-theta-formula}
	We have 
		\begin{equation} \label{eq:psi-of-theta-formula}
			\psi\circ\theta_{\pi,n} = \theta_{\pi,n}^q + \pi\theta_{\pi,n+1} 
				+ \pi P(\theta_{\pi,0},\dots,\theta_{\pi,n-1}),
		\end{equation}
	for some polynomial $P(\theta_{\pi,0},\dots,\theta_{\pi,n-1})$ with coefficients in $R$.
\end{lemma}
\begin{proof}
	It is clear for $n=0$.  For $n\geq 1$, we will use induction.  Recall the general implication
		\begin{equation*} \label{eq:witt-components-2}
		x\equiv y \bmod \m \quad\Longrightarrow\quad x^{q^{j}}\equiv y^{q^{j}} \bmod \m^{j+1},
		\end{equation*}
	for $j\geq 1$, which itself is easily proved by induction. Together with the formula 
	(\ref{eq:psi-of-theta-formula})	for $\psi\circ\theta_{\pi,i}$ with $i<n$, this implies
		\begin{align*}
		\psi\circ\psi^{\circ n} 
			&= \sum_{i=0}^n \pi^i(\psi\circ\theta_{i})^{q^{n-i}} \\
			&\equiv \pi^n\psi\circ\theta_{n}+ \sum_{i=0}^{n-1} \pi^i(\theta_{i}^q)^{q^{n-i}} 
				\mod \m^{n+1}R[\theta_{0},\dots,\theta_{n-1}]
		\end{align*}
	When this is combined with the defining formula (\ref{eq:witt-comp}) for $\psi^{\circ(n+1)}$,
	we have
		$$
			\pi^n\psi\circ\theta_{n} \equiv \pi^n\theta_{n}^q + \pi^{n+1}\theta_{n+1}
				\bmod\m^{n+1}R[\theta_{0},\dots,\theta_{n-1}].
		$$
	Dividing by $\pi^n$ completes the proof.
\end{proof}

\begin{prop} \label{pro:witt-components} 
	The elements $\theta_{\pi,0},\theta_{\pi,1},\dots$ of $R[1/\pi]\tn_R \Lambda_{\m}$ 
	lie in $\Lambda_{\m}$, and they generate $\Lambda_{\m}$ freely as an $R$-algebra.
	Further, the elements $\theta_{\pi,0},\dots,\theta_{\pi,n}$ 
	lie in $\Lambda_{\m,n}$, and they generate $\Lambda_{\m,n}$ freely as an $R$-algebra.
\end{prop}

This is essentially Witt's theorem 1~\cite{Witt:Vectors}. 

\begin{proof}
 	By induction, the elements $\theta_0,\dots,\theta_n$ generate the same sub-$R[1/\pi]$-algebra
	of $R[1/\pi]\tn_R \Lambda_{\m}$ as $\psi^{\circ 0},\dots,\psi^{\circ n}$,
	and are hence algebraically independent over $R[1/\pi]$.
	Since $R\subseteq R[1/\pi]$, they are also algebraically independent over $R$.

	Let $B_n$ be the sub-$R$-algebra of $R[1/\pi]\tn_R \Lambda_{\m}$
	generated by $\theta_{0},\dots,\theta_{n}$, and let $B=\bigcup_n B_n$.	
	To show $\Lambda_{\m}\supseteq B$, we may
	assume by induction that $\Lambda_{\m}\supseteq B_n$
	and then show $\Lambda_{\m}\supseteq B_{n+1}$.
	By \ref{lem:psi-of-theta-formula} and because $\Lambda_{\m}$ is 
	a $\Lambda_{\m}$-ring, we have 
		$$
		\pi\theta_{n+1} \in
			\big(\psi\circ\theta_{n} - \theta_{n}^q \big) + \m \Lambda_{\m,n} 
			\subseteq \m\Lambda_{\m}.
		$$
	Dividing by $\pi$, we have $\theta_{n+1}\in \Lambda_{\m}$,
	and hence $\Lambda_{\m}\supseteq B_n[\theta_{n+1}]=B_{n+1}$.

	On the other hand, by \ref{lem:psi-of-theta-formula} again, we have
		$$
		\psi\circ\theta_{n}\equiv\theta_{n}^q\bmod \m B_{n+1}
		$$
	for all $n$.
	Therefore $B$, since it is generated by the $\theta_{n}$,
	is a sub-$\Lambda_{\m}$-ring of $R[1/\pi]\tn_R \Lambda_{\m}$.
	It follows that $B\supseteq \Lambda_{\m}\nc e=\Lambda_{\m}$
	and hence that $B=\Lambda_{\m}$.  
	
	Last, the equality $\Lambda_{\m,n}=B_n$ follows immediately from the above:
		\begin{align*}
			\Lambda_{\m,n} &= \Lambda_{\m}\cap \big(R[1/\pi]\tn_R\Psi_{\m,n}\big) 
					= B \cap \big(R[1/\pi]\tn_R\Psi_{\m,n}\big) \\
				&= R[\theta_{0},\dots] \cap R[1/\pi][\theta_{0},\dots,\theta_{n}] \\
				&= R[\theta_{0},\dots,\theta_{n}] = B_n.
		\end{align*}
\end{proof}

\subsection{} \emph{Example: Presentations of $\Lambda_{\m,n}\bcp A$.}
\label{subsec:wnls-affine-presentation}
Using~\ref{pro:witt-components}, we
can turn a presentation of an $R$-algebra $A$ into a presentation of
$\Lambda_{\m,n}\bcp A$.  
We have
	$$
	\Lambda_{\m,n}\bcp R[x] \cong \Lambda_{\m,n} = R[\theta_0,\dots,\theta_n],
	$$
where $\theta_k$ is short for $\theta_{\pi,k}$, which corresponds
to the element $\theta_{\pi,k}(x)=\theta_{\pi,k}\bcp x$.

Because the functor $\Lambda_{\m,n}\bcp\vbl$ preserves coproducts and 
coequalizers, we have
	\begin{equation} \label{eq:wnls-presentation-theta}
	\Lambda_{\m,n}\bcp \big(R[x_1,\dots,x_r]/(f_1,\dots,f_s)\big )
		= R[\theta_i(x_j)]/(\theta_i(f_k)),
	\end{equation}
where $0\leq i\leq n$, $1\leq j \leq r$, and $1\leq k\leq s$.
Here each expression $\theta_i(x_j)$ is a single free variable,
and $\theta_i(f_k)$ is understood to be the polynomial in the variables
$\theta_i(x_j)$ that results from expanding $\theta_i(f_k)$ using the
sum and product laws for $\theta_i$.  Because $\Lambda_{\m,n}\bcp\vbl$ preserves
filtered colimits, we can give a similar presentation of $\Lambda_{\m,n}\bcp A$ 
for any $R$-algebra $A$.  Similarly, we can take the colimit over $n$ to get a presentation
for $\Lambda_{\m}\bcp A$.

In the $\ptst$-typical case, where $\ptst$ is finite, one can write down a presentation
of $\Lambda_{R,\ptst}\bcp A$
by iterating~(\ref{eq:wnls-presentation-theta}), according 
to~\ref{pro:lambda-iterate-generalized} below.  We can pass from the case where $\ptst$ is
finite to the case where it is arbitrary by taking colimits, 
as in~\ref{pro:lambda-for-colimit-of-ptst}.

The method above is not particular to the $\theta$ operators---it works for any subset of
$\Lambda_{\m,n}$ that generates it freely 
as an $R$-algebra. 
For example, we can use the $\delta$ operators
of~\ref{rmk:formula-based-approaches}. Let $\delta^{i}\in\Lambda_{\m}$ denote the $i$-th iterate of
$\delta_{\pi}$. Then the elements $\delta^{0},\dots,\delta^{n}$ lie in $\Lambda_{\m,n}$ and freely
generate it as an $R$-algebra. (As in~\ref{pro:witt-components}, this follows by
induction, but in this case, there are no subtle congruences to check.) Therefore we have
	\begin{equation} \label{eq:wnls-presentation-delta}
	\Lambda_{\m,n}\bcp \big(R[x_1,\dots,x_r]/(f_1,\dots,f_s)\big)
		= R[\delta^{i}(x_j)]/(\delta^{i}(f_k)),
	\end{equation}
where $0\leq i\leq n$, $1\leq j \leq r$, and $1\leq k\leq s$.
We interpret the expressions $\delta^{i}(x_j)$ and $\delta^{i}(f_k)$
as above.  
The general $\ptst$-typical case can be handled as above.
(See Buium--Simanca~\cite{Buium-Simanca:Arithmetic-Laplacians}, proof of proposition 2.12.)

\subsection{} \label{subsec:Witt-components}
{\em Witt components.} 
It follows from~\ref{pro:witt-components} that, given $\pi$, 
we have a bijection 
	\begin{equation} \label{map:witt-component-infinite-length}
		W_{\m}(A)\longisomap A\times A\times\cdots,
	\end{equation}
which sends a map $f\:\Lambda_{\m}\to A$ to the sequence 
$(f(\theta_{\pi,0}),f(\theta_{\pi,1}),\dots)$.  To make the dependence on $\pi$ explicit, we 
will often write $(x_0,x_1,\dots)_\pi$ 
for the image of $(x_0,x_1,\dots)$ under the inverse of this map.  
If $R=\bZ$ and $\pi=p$, then this identifies $W_{\m}(A)$ with the  
ring of $p$-typical Witt vectors as defined traditionally.
Similarly, when $R$ is a complete discrete valuation ring, we get an identification of
$W_{\m}(A)$ with Hazewinkel's ring of ramified Witt vectors $W^R_{q,\infty}(A)$.
(See~\cite{Hazewinkel:book}, (18.6.13), (25.3.17), and (25.3.26)(i).)
We call the $x_i$ the \emph{Witt components} (relative to $\pi$)
of the element $(x_0,\dots)_{\pi}\in W(A)$.

Similarly, using the free
generating set $\theta_{\pi,0},\dots,\theta_{\pi,n}$ of $\Lambda_{\m,n}$,
we have a bijection
	\begin{equation} \label{map:witt-component-finite-length}
		W_{\m,n}(A) \longisomap A^{[0,n]}.
	\end{equation}
As above, we will write $(x_0,\dots,x_n)_{\pi}$ for the image of $(x_0,\dots,x_n)$
under the inverse of this map.  This identifies $W_{\m,n}(A)$ with the traditionally defined
ring of $p$-typical Witt vectors of length $n+1$. 
(For remarks on the $+1$ shift, see~\ref{subsec:traditional-versus-normalized-indexing}.)

Note that the Witt components do not depend on the choice of $\pi$ in a simple, multilinear way.
For example, if $u$ is an invertible element of $R$ and we have
	$$
	(x_0,x_1,\dots)_\pi=(y_0,y_1,\dots)_{u\pi},
	$$ 
then we have
	$$
	x_0=y_0, \quad\quad
	x_1 = u y_1, \quad\quad
	x_2 = u^2 y_2 + \pi^{-1}(u-u^q)y_1^q,\quad \dots.
	$$

As in~\ref{subsec:wnls-affine-presentation}, we could use the free generating set
$\delta^0,\delta^1,\dots$ of $\Lambda_{\m}$ instead of $\theta_0,\theta_1,\dots$. This would give a
different bijection between $W_{\m}(A)$ and the set $A^{\bN}$, and hence an $R$-algebra structure
on the set $A^{\bN}$ which is isomorphic to Witt's but not equal to it. The truncated versions
agree up to $A\times A$, but differ after that. This is simply because $\delta^0=\theta_{0}$ and
$\delta^1=\theta_1$, but $\delta^2\neq \theta_2$. (See Joyal~\cite{Joyal:Witt}, p.\ 179.)

\subsection{} \emph{The ghost principle.}
It follows from the descriptions~(\ref{map:witt-component-infinite-length})
and~(\ref{map:witt-component-finite-length}) that $W_{\m}$ and $W_{\m,n}$ preserve surjectivity. On
the other hand, every $R$-algebra is a quotient of an $\m$-flat $R$-algebra (even a free one).
Therefore to prove any functorial identity involving rings of Witt vectors when $\m$ is principal,
it is enough to restrict to the $\m$-flat case. Further, any $\m$-flat $R$-algebra $A$ is contained
in an $R[1/\m]$-algebra, such as $R[1/\m]\tn_R A$. Since $W_{\m}$ and $W_{\m,n}$, being
representable functors, preserve injectivity, it is even enough to check functorial identities on
$R[1/\m]$-algebras $A$, in which case rings of Witt vectors agree with the much more tractable
rings of ghost components. An example with details is given in
\ref{subsec:local-Verschiebung-principal}.

\subsection{} {\em Verschiebung.} 
\label{subsec:local-Verschiebung-principal}
For any $R$-algebra $A$ 
define an operator $V_{\pi}$, called the Verschiebung (relative to $\pi$), 
on $W_{\m}(A)$ by 
\begin{equation} \label{eq:V-witt-component-description-new}
	V_\pi\big((y_0,y_1,\dots)_\pi\big)=(0,y_0,y_1,\dots)_\pi.
\end{equation}
This is clearly functorial in $A$. Define another, identically denoted operator
on the ghost ring $A^{\bN}$ by the formula
\begin{equation} \label{eq:V-ghost-description-new}
	V_\pi\big(\AngBrak{z_0,z_1,\dots}\big)=\AngBrak{0,\pi z_0,\pi z_1,\dots}.
\end{equation}
These operators are compatible in that
we have $\gh{}(V_{\pi}(y))=V_{\pi}(\gh{}(y))$ for all $y\in W_{\m}(A)$, and
the operator $V_{\pi}$ on the ghost ring is clearly $R$-linear.
It follows by the ghost principle that the operator
$V_{\pi}$ on $W_{\m}(A)$ is $R$-linear. Here is the argument in some detail.

We need to check the identities $rV_{\pi}(y)=V_{\pi}(ry)$
and $V_{\pi}(x+y)=V_{\pi}(x)+V_{\pi}(y)$, for $r\in R$, $x,y\in W_{\m}(A)$.
Write $x=(x_0,x_1,\dots)_{\pi}$ and $y=(y_0,y_1,\dots)_{\pi}$.
If $A$ is a $\ptst$-flat, the ghost map $\gh{}\:W_{\m}(A)\to A^{\bN}$ is injective.
Therefore $V_{\pi}$ is $R$-linear on $W_{\m}(A)$, by the $R$-linearity of $V_{\pi}$ on 
the ghost ring.

The general case then follows from $\ptst$-flat case.
Fix an $\ptst$-flat $R$-algebra $\tilde{A}$ with a surjective $R$-algebra map
$\tilde{A}\to A$. For each $i$, let $\tilde{y}_i$ 
be a pre-image of $y_i$,
and set $\tilde{y}=(\tilde{y}_0,\dots)_{\pi}\in W_{\m}(A)$.
The induced
map $f\:W_{\m}(\tilde{A})\to W_{\m}(A)$ then satisfies $f(\tilde{y})=y$. Therefore we have
	\begin{align*}
	V_{\pi}(ry) &= V_{\pi}(rf(\tilde{y})) = V_{\pi}(f(r\tilde{y})) = f(V_{\pi}(r\tilde{y})) \\
		&= f(rV_{\pi}(\tilde{y})) = rf(V_{\pi}(\tilde{y})) = rV_{\pi}(f(\tilde{y})) = rV_{\pi}(y).
	\end{align*}
The additivity axiom follows similarly.

\subsection{} \emph{Example.}
\label{subsec:W_n-of-R-example}
$W_{R,\m,n}(R)$ has a presentation
	\[
	R[x_1,\dots,x_{n}]/(x_ix_j-\pi^ix_j\,|\,1\leq i\leq j\leq n),
	\]
where the element $x_i$ corresponds to $V_{\pi}^i(1)$.

\subsection{} {\em Teichm\"uller lifts.}
\label{subsec:teich-2}
Under the composition
	$$
	\xymatrix@C=30pt{
	A \ar^-{a\mapsto[a]}[r]
		& W(A) \ar^-{\gh{}} [r]
		& A \times A \times \cdots
	}
	$$
(see~\ref{subsec:teich-1}),
the image of $a$ is $\AngBrak{a,a^q,a^{q^2},\dots}$.  It follows from the ghost principle that
	\[
	[a] = (a,0,0,\dots)_{\pi} \in W(A).
	\]
Multiplication by Teichm\"uller lifts also has a simple description in terms of
Witt components:
	\begin{equation} \label{eq:teich-mult-on-witt-comp}
		[a] (\dots,b_i,\dots)_{\pi} = (\dots,a^{q^i}b_i,\dots)_{\pi}.
	\end{equation}
Again, this follows from the ghost principle.

\section{General single-prime case}
\label{sec:general-single-prime-case}

Assume $\ptst$ consists of a single ideal $\m$, possibly not principal. 
Let $n$ be an element of $\bN$. Let us write $W_{R,\m,n}=W_{R,\ptst,n}$ and so on.

Let $K_{\m}$ denote $R_{\m}[1/\m]$. If $\m$ is the unit ideal, we understand $R_{\m}$, and hence
$K_{\m}$, to be the zero ring. Otherwise, $R_{\m}$ is a discrete valuation ring and $K_{\m}$ is its
fraction field. In particular, $\m$ becomes principal in $R[1/\m]$, $R_{\m}$, and $K_{\m}$.
The following proposition then allows us to describe $W_{R,\m,n}(A)$ 
in terms of the case where $\m$ is principal, and hence in terms of Witt components.

\begin{proposition}\label{pro:principal-principle}
	For $R'=R[1/\m],R_{\m},K_{\m}$, write $W_{R',\m,n}=W_{R',\m R',n}$.
	Then for any $R$-algebra $A$, the ring $W_{R,\m,n}(A)$ is the equalizer of the two maps
		$$
			\displayrightarrows{W_{R[1/\m],\m,n}\big(R[1/\m]\tn_R A\big)\times
				W_{R_{\m},\m,n}\big(R_{\m}\tn_R A\big)}
			{W_{K_{\m},\m,n}(K_{\m}\tn_R A)}		
		$$
	induced by projection onto the two factors and the bifunctoriality of $W_{\vbl,\m,n}(\vbl)$.
\end{proposition}

\begin{proof}
The diagram
	$$
	\displaylabelfork{R}{}{R[1/\m]\times R_{\m}}{\pr_1}{\pr_2}{K_{\m}}
	$$
is an equalizer diagram.  Since $K_{\m}$ is $\m$-flat, so is any sub-$R$-module of $K_{\m}$.
It follows that for any $R$-algebra $A$, the induced diagram
	$$
	\displayfork{A}{\Big(R[1/\m]\times R_{\m}\Big)\tn_R A}{K_{\m}\tn_R A}
	$$
is an equalizer diagram.  Since 
$W_{R,\m,n}$ is representable, it preserves equalizers,
and so the induced diagram (writing $W_n=W_{R,\m,n}$)
	\begin{equation*} 
		\displayfork{W_{n}(A)}
			{W_{n}\big(R[1/\m]\tn_R A\big)\times
				W_{n}\big(R_{\m}\tn_R A\big)}
			{W_{n}(K_{\m}\tn_R A)}		
	\end{equation*}
is also an equalizer diagram.
Then (\ref{map:witt-independence-of-R-finite-length}) completes the proof.
\end{proof}

\subsection{} \emph{Verschiebung in general.}
We can define Verschiebung maps
	\begin{equation} \label{map:naturalized-Verschiebung}
		V^j\:\m^j\tn_R W_{R,\m}(A) \longmap W_{R,\m}(A).
	\end{equation}
To do this,
it is enough, by~\ref{pro:principal-principle}, to restrict to the case where $\m$
is principal, as long as our construction is functorial in $A$ and $R$.
So, choose a generator $\pi\in\m$ and define
	\begin{equation}
		V^j(\pi^j\tn y)= V_{\pi}^j(y),
	\end{equation}
for all $y\in W_{R,\m}(A)$.
On ghost components it satisfies 
	$$
	V^j(x\tn\AngBrak{z_0,z_1,\dots}) =\AngBrak{0,\dots,0,xz_0,xz_1,\dots},
	$$
where the number of leading zeros is $j$.  In particular, it is independent
of the choice of $\pi$, by the ghost principle.

If we write $W_{R,\m}(A)_{(j)}$ for $W_{R,\m}(A)$ viewed as a $W_{R,\m}(A)$-algebra by way
of the map $\psi_j\:W_{R,\m}(A)\to W_{R,\m}(A)$, then the map
	\begin{equation}
		V^j\: \m^j\tn_R W_{R,\m}(A)_{(j)} \longmap W_{R,\m}(A),
	\end{equation}
is $W_{R,\m}(A)$-linear, as is easily checked using the ghost principle.
Expressed as a formula, it says
	\begin{equation}
		\quad V^j(x\tn y\psi_j(z)) = V^j(x\tn y)z.		
	\end{equation}
In particular, the image $V^jW_{R,\m}(A)$ of $V^j$ is an ideal of $W_{R,\m}(A)$.

Let us also record the identities
	\begin{equation} \label{eq:FV=p}
		\psi_j\big(V^j(x\tn y)\big) = xy
	\end{equation}
and
	\begin{equation} \label{eq:Verschiebung-product-new}
		V^j(x\tn y)V^{j}(x'\tn y')= xV^{j}(x'\tn yy') \in\m^j V^jW_{R,\m}(A).
	\end{equation}
Again, one checks these using the ghost principle. 

Finally, for any $n\in\bN$, the map $V^j$ descends to a map
	\begin{equation} \label{def:Verschiebung-finite-length}
		V^j\:\m^j\tn_R W_{R,\m,n}(A)_{(j)} \longmap W_{R,\m,n+j}(A),
	\end{equation}
and the obvious analogues of the identities above hold here.

\subsection{} \emph{Remark.}
We can define Verschiebung maps even if we no longer assume there is only one ideal in $\ptst$.
For any $j\in\nset$, let $J$ denote the ideal $\prod_{\pind}\m_{\pind}^{j_{\pind}}$ of $R$.
Then $V^j$ would be a map $J\tn_R W_{R,\ptst}(A)\to W_{R,\ptst}(A)$.  The identities above, 
suitably interpreted, continue to hold.  We will not need this multiple-prime version.

\begin{proposition}\label{pro:V-exact-sequence}
	The sequence
		\begin{equation} \label{diag:V-exact-sequence}
			0 \longmap \m^j\tn_R W_{R,\m,n}(A)_{(j)} \longlabelmap{V^j} W_{R,\m,n+j}(A) 
				\longmap W_{R,\m,j}(A) \longmap 0
		\end{equation}
	is exact.
\end{proposition}
\begin{proof}
	Write $W_{R',n}=W_{R',\m R',n}$ when $R'$ is an $R$-algebra such
	that the ideal $\m R'$ is supramaximal.
	
	First consider the case where $\m$ is principal.  Let $\pi\in\m$ be a generator.
	Using~(\ref{eq:V-witt-component-description-new}), it is clear
	that $V^j$ is injective and that its image is the set of
	Witt vectors whose Witt components (relative $\pi$) are $0$ in 
	positions $0$ to $j-1$.  By~\ref{subsec:Witt-components}, 
	the pre-image of $0$ under the map $W_{R,n+j}(A)\to W_{R,j}(A)$ is the same subset, and
	the map $W_{R,n+j}(A)\to W_{R,j}(A)$ is surjective.
	
	Now consider the general case.  Augment the diagram~(\ref{diag:V-exact-sequence})
	by expressing each term of~(\ref{diag:V-exact-sequence})
	as an equalizer as in~\ref{pro:principal-principle}.
	Here we use that $\m$ is $R$-flat.  It then follows from the principal
	case and the snake lemma that~(\ref{diag:V-exact-sequence}) is left exact.
	
	It remains to prove that the map $W_{R,n+j}(A)\to W_{R,j}(A)$ is surjective.
	By induction, we can assume $n=1$.  
	By~\ref{pro:principal-principle}, for any $i\in\bN$ we have
		\[
		W_{R,i}(A)= W_{R_\m,i}(R_\m\tn_R A) \times_{W_{K_\m,i}(K_{\m}\tn_R A)} 
			W_{R[1/\m],i}(R[1/\m]\tn_R A).
		\]
	Now let $\pi$ denote a generator of the maximal ideal of $R_{\m}$, and
	suppose two elements
	   \begin{align*}
	   y &= (y_0,\dots,y_{j})_\pi\in W_{R_\m,j}(R_\m\tn_R A),\\
	   z &= \AngBrak{z_0,\dots,z_{j}}\in (R[1/\m]\tn_R A)^{j+1}=W_{R[1/\m],j}(R[1/\m]\tn_R A)
	   \end{align*}
	have the same image in $W_{K_\m,j}(K_{\m}\tn_R A)$.  To lift the corresponding element of 
	$W_{j}(A)$ to $W_{j+1}(A)$, we need to find elements 
		$$
		y_{j+1}\in R_\m\tn_R A\quad \text{and}\quad z_{j+1}\in R[1/\m]\tn_R A
		$$ 
	such that in $K_{\m}\tn_R A$, we have
		\begin{equation}
		\label{eq:witt-lift}
		y_0^{q^{j+1}}+\cdots+\pi^{j+1} y_{j+1} = z_{j+1}.
		\end{equation}
	So, choose an element $z_{j+1}\in A$ whose image under the surjection
		$$
		A \longmap A/(\m A)^{j+1}= R_\m/(\m R_{\m})^{j+1}\tn_R A
		$$ 
	agrees with the image of $y_0^{q^{j+1}}+\cdots+\pi^{j}y_{j}$.
	It follows that the element
		\[
		y_0^{q^{j+1}}+\cdots+\pi^{j}y_{j}^q-1\tn z_{j+1} \in R_\m\tn_R A
		\]
	lies in $\pi^{j+1}(R_\m\tn_R A)$.  It thus equals $\pi^{j+1}y_{j+1}$ for some
	element $y_{j+1}\in R_\m\tn_R A$.
	And so $y_{j+1}$ and $z_{j+1}$ satisfy~(\ref{eq:witt-lift}).
\end{proof}

\begin{corollary}\label{cor:gr-V}
	For any $R$-algebra $A$, we have
		\begin{equation}
			\bigoplus_{i\in[0,n]} \m^i\tn_R A_{(i)} \longisomap \gr_V W_{R,\m,n}(A),
		\end{equation}
	where $A_{(i)}$ denotes $A$ viewed as a $W_n(A)$-module via the
	ring map $\gh{i}\:W_n(A)\to A$.
\end{corollary}

\subsection{} {\em Reduced ghost components.}
\label{subsec:generalized-ghost-definition}
We can define infinitely many ghost components for Witt vectors of finite
length $n$ if we are willing to settle for answers modulo $\m^{n+1}$.

First assume $\m$ is generated by some element $\pi$.  By examining the Witt
polynomials~(\ref{eq:witt-comp}), we can see that for any $i\geq 0$, the composition
\[
W_{R,\m}(A) \longlabelmap{\gh{i}} A \longmap A/\m^{n+1}A
\]
vanishes on $V^{n+1}W_{R,\m}(A)$.  It therefore
factors through $W_{R,\m,n}(A)$, giving a map $\rgh{i}$ from $W_{R,\m,n}(A)$ to $A/\m^{n+1}A$.  

When $\m$ is not assumed to be principal, we define $\rgh{i}$ by
localizing at $\m$:
	$$
	W_{R,\m,n}(A) \longmap W_{R_{\m},\m R_{\m},n}(R_{\m}\tn_R A) \longlabelmap{\rgh{i}} 
		(R_{\m}\tn_R A)/\m^{n+1}(R_{\m}\tn_R A) = A/\m^{n+1}A,
	$$
where the middle map is $\rgh{i}$ as constructed above in the principal case.
We call the composition
	\begin{equation} \label{map:affine-reduced-ghost-map}
		W_{R,\m,n}(A)\longlabelmap{\rgh{i}} A/\m^{n+1}A
	\end{equation}
the $i$-th {\em reduced ghost component} map.

\section{Multiple-prime case}
\label{sec:multiple-prime-case}

The purpose of this section is to give some results on reducing the family $\ptst$
(of~\ref{subsec:affine-general-notation}) to simpler families. The first reduces from the case
where $\ptst$ is arbitrary to the case where it is finite, and the second reduces from the case
where it is finite to the case where it has a single element. We will often write 
$W_{\ptst}=W_{R,\ptst}$, $\Lambda_{\ptst}=\Lambda_{R,\ptst}$, and so on, for short.

\begin{proposition} \label{pro:lambda-for-colimit-of-ptst}
	The canonical maps
	\begin{align}
			\colimm_{\ptst'} \Lambda_{R,\ptst'} &\longmap \Lambda_{R,\ptst},
			 	\label{map:lambda-for-colimit-of-ptst-infinite}\\
			\colimm_{\ptst'} \Lambda_{R,\ptst',n'} &\longmap \Lambda_{R,\ptst,n},
				\label{map:lambda-for-colimit-of-ptst-finite}
	\end{align}
	are isomorphisms.  Here $\ptst'$ runs over the finite subfamilies of $\ptst$,
	and $n'$ is the restriction to $\ptst'$ of a given element $n\in\nset$.
\end{proposition}
\begin{proof}
	Consider (\ref{map:lambda-for-colimit-of-ptst-infinite}) first.
	Since each map $\Lambda_{\ptst'}\to\Lambda_{\ptst}$ is an injection,
	(\ref{map:lambda-for-colimit-of-ptst-infinite}) is an injection.
	Therefore, since $\Lambda_{\ptst}$ is freely generated as a $\Lambda_{\ptst}$-ring by
	the element $e=\psi_0$, it is enough to show the sub-$\Psi_{\ptst}$-ring
	$\colimm_{\ptst'}\Lambda_{\ptst'}$ of $\Lambda_{\ptst}$ 
	is a sub-$\Lambda_{\ptst}$-ring.
	Since it is flat, we only need to check the Frobenius lift property.
	So, suppose $\m\in\ptst$.  For any element $x$ of the colimit, there is a
	finite family $\ptst''$ such that $x\in\Lambda_{\ptst''}$ and $\m\in\ptst''$.
	But $\Lambda_{\ptst''}$ is a $\Lambda_{\ptst''}$-ring. So we have
	$\psi_{\m}(x)\equiv x^{q_{\m}}$ modulo $\m\Lambda_{\ptst''}$, and hence
	modulo $\m(\colimm_{\ptst'}\Lambda_{\ptst'})$. Therefore the Frobenius lift
	property holds for the colimit ring.

	Then (\ref{map:lambda-for-colimit-of-ptst-finite}) follows:
	\begin{align*}
		\Lambda_{\ptst,n} &= (R[1/\ptst]\tn_R \Psi_{\ptst,n}) \cap \Lambda_{\ptst} 
			= (\colimm_{\ptst'} R[1/\ptst]\tn_R \Psi_{\ptst',n'})\cap 
				\colimm_{\ptst'}\Lambda_{\ptst'} \\
			&= \colimm_{\ptst'} \big((R[1/\ptst]\tn_R \Psi_{\ptst',n'})\cap \Lambda_{\ptst'}\big)
			= \colimm_{\ptst'} \Lambda_{\ptst',n'}.
	\end{align*}
\end{proof}

\begin{corollary} \label{cor:W-for-colimit-of-ptst}
	For any $R$-algebra $A$, the canonical maps
		\begin{align}
				\label{map:W-for-colimit-of-ptst-1}
			W_{R,\ptst}(A) &\longmap \lim_{\ptst'} W_{R,\ptst'}(A), \\
				\label{map:W-for-colimit-of-ptst-2}
			W_{R,\ptst,n}(A) &\longmap \lim_{\ptst'} W_{R,\ptst',n'}(A)			
		\end{align}
	are isomorphisms, where $\ptst'$, $n$, and $n'$ are as in 
	\ref{pro:lambda-for-colimit-of-ptst}.
\end{corollary}

\begin{proposition}\label{pro:lambda-iterate-generalized}
Let $\ptst'\smcoprod\ptst''$ be a partition of $\ptst$.
Then the canonical maps
	\begin{align}  
		\Lambda_{R,\ptst'}\bcp_R\Lambda_{R,\ptst''} &\longmap \Lambda_{R,\ptst} 
			\label{map:lambda-iterate-generalized-infinite} \\
		\Lambda_{R,\ptst',n'}\bcp_R\Lambda_{R,\ptst'',n''} &\longmap \Lambda_{R,\ptst,n} 
			\label{map:lambda-iterate-generalized-finite} 
	\end{align}
are isomorphisms, where $n'$ and $n''$ denote the restrictions to $\ptst'$ and $\ptst''$
of a given element $n\in\nset$.
\end{proposition}
\begin{proof}
	It is enough to show each map becomes an isomorphism after base change to $R[1/\ptst']$
	and $R[1/\ptst'']$. So, by~(\ref{map:lambda-change-of-R-infinite-length}),
	we can assume every element in either $\ptst'$ or $\ptst''$ is the unit ideal.
	
	In the second case, we have
		$$
		\Lambda_{\ptst'} \bcp_R \Lambda_{\ptst''} 
			= \Lambda_{\ptst'} \bcp_R R[\bN^{(\ptst'')}]
			= \Lambda_{\ptst'}[\bN^{(\ptst'')}]
			= \Lambda_{\ptst}
		$$
	The argument for (\ref{map:lambda-iterate-generalized-finite}) is the same,
	but we replace the generating set $\bN^{(\ptst'')}$ with $[0,n'']$.
	
	Now suppose every element in $\ptst'$ is the unit ideal.
	Then a $\Lambda_{\ptst'}$-ring is the same as a $\Psi_{\ptst'}$-ring. So
	we have
		$$
		\Lambda_{\ptst'}\bcp_R\Lambda_{\ptst''} = \Lambda_{\ptst''}[\bN^{(\ptst')}]
			= \Lambda_{\ptst}.
		$$
	For (\ref{map:lambda-iterate-generalized-finite}), replace $\bN^{(\ptst')}$
	with $[0,n']$, as above.
\end{proof}

\begin{corollary}\label{cor:witt-iterate-generalized}
Let $\ptst'\smcoprod\ptst''$ be a partition of $\ptst$.
Then for any $R$-algebra $A$, the canonical maps
	\begin{align}
			\label{map:witt-iterate-generalized-infinite-length}
		W_{R,\ptst}(A) &\longmap W_{R,\ptst''}\big(W_{R,\ptst'}(A)\big) \\
			\label{map:witt-iterate-generalized-finite-length}
		W_{R,\ptst,n}(A) &\longmap W_{R,\ptst'',n''}\big(W_{R,\ptst',n'}(A)\big)		
	\end{align}
are isomorphisms, where $n$, $n'$, $n''$ are as in \ref{pro:lambda-iterate-generalized}.	
\end{corollary}

\subsection{} \emph{Remark.}
By the results above, it is safe to say that expressions such as
	\begin{equation} \label{eq:independent-of-ordering}
	\Lambda_{\m_1}\bcp_R\cdots\bcp_R\Lambda_{\m_r} \quad\text{ and } \quad
		W_{\m_r}\circ\cdots\circ W_{\m_1}(A)		
	\end{equation}
are independent of the ordering of the $\m_i$, assuming the $\m_i$ are pairwise coprime. 
(Note that it is not generally true
that $P\bcp P' \cong P'\bcp P$ for plethories $P$ and $P'$. See~\cite{Borger-Wieland:PA}, 2.8.)

If we ask that the expressions in (\ref{eq:independent-of-ordering}) be independent only up to
isomorphism, then it is not even necessary that the $\m_{\pind}\in\ptst$ be pairwise coprime
(\ref{subsec:affine-general-notation}). But invariance up to isomorphism is not a such a useful
property, and most of the time coprimality really is necessary. For example, we could look at
rings with more than one Frobenius lift at a single maximal ideal, but we would not be able to
reduce to the case of a single Frobenius lift. Indeed, if $\ptst$ consists of a single maximal
ideal $\m$, the two endomorphisms $\psi_{WW(A)}$ and $W(\psi_{W(A)})$ of $WW(A)$ commute, and the
first is clearly a Frobenius lift, but the second is generally not. Therefore $WW(A)$ cannot be the
cofree ring with two commuting Frobenius lifts at $\m$.

In fact, I believe this is the only place where we use the coprimality assumption directly. The
rest of our results depend on it only through~\ref{pro:lambda-iterate-generalized}. Although I know
of no applications, it would be interesting to know whether the abstract set up of this paper, and
then the main results, hold when we allow more than one Frobenius lift at each maximal ideal.

\section{Basic affine properties}
\label{sec:basic-affine-properties}

This section provides some basic results about the commutative algebra of Witt vectors.
They are just the ones needed to be able to prove the main theorems
in sections \ref{sec:affine-ghost-descent-single-prime-case} and
\ref{sec:W-and-etale-morphisms} and to set up the global theory in \cite{Borger:BGWV-II}.
There are other basic results that could have been included here, but which I have put off 
to~\cite{Borger:BGWV-II}, where they will be proved for all algebraic spaces.

We continue with the notation of~\ref{subsec:affine-general-notation}.
Fix an element $n\in\nset$. We will often write $W_n=W_{\ptst,n}=W_{R,\ptst,n}$
and so on, for short.
By~\ref{cor:W-for-colimit-of-ptst}, we may assume that $\ptst$ agrees with the
support of $n$, and in particular that it is finite.

\begin{proposition}\label{pro:Zariski-localization-on-R}
	Let $R'$ be an $\ptst$-flat $R$-algebra such that the structure map $R\to R'$ is a
	ring epimorphism (as in~\ref{subsec:localizing-R-infinite-length}).
	Then the composition
		$$
		\xymatrix@C=30pt{
		R'\tn_R W_{R,\ptst,n}(A) \ar^-{(\ref{map:witt-change-of-R-finite-length})}[r] 
			& W_{R,\ptst,n}(R'\tn_R A)
					\ar_-{\sim}^{(\ref{map:witt-independence-of-R-finite-length})^{-1}\ }[r] 
			& W_{R',\ptst',n}(R'\tn_R A)
		}
		$$
	is an isomorphism, where $\ptst'$ is as in \ref{subsec:localizing-R-infinite-length}.
\end{proposition}
\begin{proof}
	We may assume by~\ref{cor:witt-iterate-generalized}
	that $\ptst$ consists of a single ideal $\m$.
	Using  \ref{pro:principal-principle} and the flatness of $R'$ over $R$,
	we are reduced to showing that the functors
	$W_{R[1/\m],\m,n}$, $W_{R_{\m},\m,n}$, and $W_{K_{\m},\m,n}$ commute with the functor 
	$R'\tn_R\vbl$.
	Therefore we may assume that the ideal $\m$ is principal. 
	
	Write $W_n=W_{R,\m,n}$.
	The result is clear for $n=0$, because $W_0$ is the identity functor.
	So assume $n\geq 1$.  By~\ref{pro:V-exact-sequence},
	we have the following map of exact sequences
		$$
		\xymatrix@C=30pt{
			0 \ar[r]
				& R'\tn_R\m\tn_R W_{n-1}(A) \ar^-{\id_{R'}\tn V^1}[r]\ar[d]
				& R'\tn_R W_n(A) \ar[r]\ar[d]
				& R'\tn_R A \ar[r]\ar@{=}[d]
				& 0 \\
			0 \ar[r]
				& \m\tn_R W_{n-1}(R'\tn_R A) \ar^-{V^1}[r]
				& W_n(R'\tn_R A) \ar[r]
				& R'\tn_R A \ar[r]
				& 0 ,
		}
		$$
	where the vertical maps are given by~(\ref{map:witt-change-of-R-finite-length}).
	By induction the leftmost vertical arrow is an isomorphism.  Therefore
	the inner one is, too.
\end{proof}

\begin{proposition}\label{pro:W-of-ideal-product}
	For any ideal $I$ in an $R$-algebra $A$, let $W_{R,\ptst,n}(I)$ denote the kernel
	of the canonical map $W_{R,\ptst,n}(A)\to W_{R,\ptst,n}(A/I)$.  Then we have
		$$
		W_{R,\ptst,n}(I)W_{R,\ptst,n}(J)\subseteq W_{R,\ptst,n}(IJ)
		$$
	for any ideals $I$, $J$ in $A$.
\end{proposition}
\begin{proof}
	Let us first show that we may assume $\ptst$ consists of a single ideal $\m$. 
	In doing this, it will be convenient to prove the following equivalent form of the statement:
	if $IJ\subseteq K$, where $K$ is an ideal in $A$, then $W_n(I)W_n(J)\subseteq W_n(K)$.
	Suppose	$\ptst=\ptst'\smcoprod \{\m\}$.  Let $n'$ be the restriction of $n$ to $\ptst$.
	Let $I'=W_{\ptst',n'}(I)$, $J'=W_{\ptst',n'}(J)$, and $K'=W_{\ptst',n'}(K)$.
	By~\ref{cor:witt-iterate-generalized}, we have 
	$W_{\ptst,n}=W_{\m,n_{\m}}\circ W_{\ptst',n'}$,
	and hence $W_{\ptst,n}(I)=W_{\m,n_{\m}}(I')$ and so on.
	By induction, we have $I'J'\subseteq K'$, and then applying the result in the
	single-ideal case gives
		$$
		W_{\ptst,n}(I)W_{\ptst,n}(J) = W_{\m,n_{\m}}(I')W_{\m,n_{\m}}(J')
			= W_{\m,n_{\m}}(K') = W_{\ptst,n}(K).
		$$
	So we will assume $\ptst=\{\m\}$ and drop $\ptst$ from the notation.
	
	By~\ref{pro:Zariski-localization-on-R}, the statement is Zariski local on $R$, and
	so we may assume the ideal $\m$ is generated by some element $\pi$.
	We will work with Witt components relative to $\pi$.
	
	We need to show that for any elements $x=(x_0,\dots,x_n)_{\pi}\in W_n(I)$
	and $y=(y_0,\dots,y_n)_{\pi}\in W_n(J)$, the product $xy$ is in $W_n(IJ)$.
	So it is sufficient to show this in the universal case,
	where $A$ is the free polynomial algebra $R[x_0,y_0,\dots,x_n,y_n]$, $I$ is the
	ideal $(x_0,\dots,x_n)$, and $J$ is the ideal $(y_0,\dots,y_n)$.
	
	Consider the following commutative diagram
		$$
		\entrymodifiers={+!!<0pt,\fontdimen22\textfont2>}
		\def\objectstyle{\displaystyle}
		\xymatrix{
			W_n(A) \ar^{\gh{\leq n}}[r]\ar[d]
				& A^{[0,n]} \ar[d] \\
			W_n(A/IJ) \ar^{\gh{\leq n}}[r]
				& (A/IJ)^{[0,n]}.
		}
		$$
	We want to show that the image of $xy$ in $W_n(A/IJ)$ is zero.
	Since $A/IJ$ is flat (even free) over $R$, the lower map $\gh{\leq n}$ is injective,
	and so it is enough to show the image of $xy$ in $(A/IJ)^{[0,n]}$ is zero.
	But by the naturality of the ghost map, we have $\gh{\leq n}(x)\in I^{[0,n]}$
	and	$\gh{\leq n}(y)\in J^{[0,n]}$. Therefore $\gh{\leq n}(xy)$ lies in $(IJ)^{[0,n]}$,
	which maps to zero in $(A/IJ)^{[0,n]}$.
\end{proof}

\subsection{} \emph{Remark.} Although the proof of~\ref{pro:W-of-ideal-product} 
given above uses some
properties specific to Witt vector functors, the result is true for any representable ring-valued
functor. See Borger--Wieland~\cite{Borger-Wieland:PA}, 5.5.

\begin{corollary}\label{cor:W-of-nilpotent-ideal}
	Let $I$ be an ideal in an $R$-algebra $A$. If $I^m=0$, then $W_{R,\ptst,n}(I)^m=0$.
\end{corollary}

\begin{proposition}\label{pro:W-preserves-surj-of-rings}
	Let $\varphi\:A\to B$ be a map of $R$-algebras. 
	If it is surjective, then so is the map $W_{R,\ptst,n}(\varphi)\:W_{R,\ptst,n}(A)\to W_{R,\ptst,n}(B)$.
\end{proposition}
\begin{proof}
	By~\ref{cor:witt-iterate-generalized}, we
	may assume $\ptst$ consists of one ideal $\m$.
	Since surjectivity can be checked Zariski locally on $R$, it is enough
	by~\ref{pro:Zariski-localization-on-R} to assume $\m$ is principal.
 	Then using the Witt components, we can identify the set map underlying
	$W_{R,\ptst,n}(\varphi)$ with the map $\varphi^{[0,n]}\:A^{[0,n]}\to B^{[0,n]}$, which is clearly
	surjective. 
\end{proof}

\begin{corollary}\label{cor:W-preserves-quotients-by-equivalence-relations}
	If $\varphi\:A\to B$ is surjective, then
		$$
		\displaylabelcofork{W_{R,\ptst,n}(A\times_B A)}
			{W_{n}(\pr_1)}{W_{n}(\pr_2)}{W_{R,\ptst,n}(A)}{W_{n}(\varphi)}{W_{R,\ptst,n}(B)}
		$$
	is a coequalizer diagram.
\end{corollary}
\begin{proof}
	The functor $W_n$ is representable, and hence commutes with
	limits.  (See~\ref{subsec:W-finite-length}.)  Therefore 
	$W_n(A\times_B A)$ agrees with $W_n(A)\times_{W_n(B)} W_n(A)$, which
	is an equivalence relation on $W_n(A)$, the quotient by which is the image
	of $W_n(\varphi)$.  
	By~\ref{pro:W-preserves-surj-of-rings},
	this is all of $W_n(B)$.
\end{proof}

\subsection{} \emph{Remark.} This result is particularly appealing when $A$ is $\ptst$-flat
and $B$ is not.  Then we can describe $W_n(B)$ in terms of 
$W_n(A)$ and $W_n(A\times_B A)$, which are directly accessible
because $A$ and $A\times_B A$ are $\ptst$-flat.

\begin{proposition}\label{pro:W-nilpotent-at-p}
	Suppose $\ptst$ consists of one ideal $\m$, and let $A$ be an $R$-algebra.  	
	For any $i\geq 0$, the map $\Spec(\id\tn\rgh{i})$ of schemes induced by the ring map
		$$
		\id\tn\rgh{i}\:R/\m \tn_R W_{R,\ptst,n}(A) \longmap R/\m\tn_R A/\m^{n+1}
		$$
	is a universal homeomorphism. For $i=0$,
	it is a closed immersion defined by a square-zero ideal.	
\end{proposition}
\begin{proof}
	Write $W_n=W_{R,\ptst,n}$ and so on.
	Consider the diagram
		$$
		\xymatrix@C=40pt{
		R/\m\tn_R W_{n}(A)\ar^-{\id\tn\rgh{i}}[r]\ar^{\id\tn\gh{0}}[d]
			& R/\m\tn_R A/\m^{n+1}A \ar^{\sim}_{r\tn a\mapsto ra}[d] \\
		A/\m A \ar^-{x\mapsto x^{q^i}}[r] 
			& A/\m A.
		}
		$$
	To show it commutes, it is enough to assume $\m$ is principal, generated by $\pi$.
	Then commutativity follows from the obvious congruence
		$$
		\gh{i}(a)=a_0^{q^i}+\pi a_1^{q^{i-1}} + \dots+\pi^i a_i \equiv a_0^{q^i} \bmod \m A,
		$$
	for any element $a=(a_0,a_1,\dots)_{\pi}\in W(A)$.

 	Therefore, $\id\tn\rgh{i}$ is the composition of a map whose kernel
	is a nil ideal and a power of the Frobenius map. The scheme maps induced
	by both of these are universal homeomorphisms.
	
	Now let us show that $\id\tn \gh{0}$ (which equals $\id\tn\rgh{0}$)
	is a surjection with square-zero kernel.
	The map $\id\tn\gh{0}$ is surjective by~\ref{subsec:teich-1} (or~\ref{pro:V-exact-sequence}).
	So let us show the square of its kernel is zero.
	By~\ref{pro:V-exact-sequence}, the kernel of the map 
	$W_n(A)\to R/\m\tn_R A$ is the ideal $V^1W_n(A)+\m W_n(A)$.  Hence it
	is enough to show $\left(V^1W_n(A)\right)^2\subseteq \m W_n(A)$.  
	This follows from~(\ref{eq:Verschiebung-product-new}).
\end{proof}

\begin{proposition}\label{pro:W-preserves-family-surj-on-spec}
	Let $(B_i)_{i\in I}$ be a family of $A$-algebras such that the induced map
	$\coprod_i \Spec B_i \to \Spec A$ is surjective. Then the induced map 
	$$
	\coprod_i \Spec W_{R,\ptst,n}(B_i) \to \Spec W_{R,\ptst,n}(A)
	$$
	is surjective.
\end{proposition}
\begin{proof}
	By \ref{cor:witt-iterate-generalized}, it is enough to assume $\ptst$ 
	consists of one ideal $\m$.
	Further, it is enough to show surjectivity after base change to $R[1/\m]$ and to $R/\m$. 
	For $R[1/\m]$, it follows from~\ref{pro:Zariski-localization-on-R} and the equality
	$W_n(C)=C^{[0,n]}$, when $\m$ is the unit ideal.  Now
	consider base change to $R/\m$. By
	\ref{pro:W-nilpotent-at-p}, the ring $W_n(A)/\m W_n(A)$ 
	is a nilpotent extension of $A/\m A$, 
	and likewise for each $B_i$, and so we are reduced to showing that
	$\coprod_i \Spec B_i/\m B_i \to \Spec A/\m A$ is surjective. 
	This is true since base change distributes over disjoint unions and preserves surjectivity.
\end{proof}

\begin{proposition}\label{pro:W-preserves-filtered-colimits}
	The $R$-algebra $\Lambda_{R,\ptst,n}$ is finitely presented,
	and the functor $W_{R,\ptst,n}$ preserves filtered colimits of $R$-algebras.
\end{proposition}
\begin{proof}
	Since $W_{R,\ptst,n}$ is represented by $\Lambda_{R,\ptst,n}$, the two statements to be proved
	are equivalent.	By~\ref{cor:witt-iterate-generalized}, we may assume $\ptst$ consists of a 
	single ideal $\m$. By EGA IV 2.7.1~\cite{EGA-no.24}, the first statement can be 
	verified fpqc locally on $R$, and in particular after base change to $R[1/\m]$ and to 
	$R_{\m}$.  Therefore by~(\ref{map:lambda-change-of-R-finite-length}), we can assume $\m$ is 
	generated by a single element $\pi$. But by~\ref{pro:witt-components},
	the $R$-algebra $\Lambda_{R,\ptst,n}$ is generated by the
	finite set $\theta_{\pi,0},\dots,\theta_{\pi,n}$.
\end{proof}

\section{Some general descent}
\label{sec:some-general-descent}

The purpose of this section is to record some facts about descent of \'etale algebras which we will
use to prove the main theorem (\ref{thm:witt-etale-general}). The results 
mention nothing about Witt vectors or anything else in this paper. So it would
be reasonable to skip this section and refer back to it only as needed.

More precisely, we do the following. First, we set up some language and notation for descent. This
is essentially a repetition of parts of Grothendieck's TDTE I~\cite{TDTE-I}. (It could not be
anything but.) Second, we prove an abstract result (\ref{pro:abstract-descent=gluing}) relating
gluing data and descent data for certain simple gluing constructions. Third, we recall
Grothendieck's theorem (\ref{thm:Grothendieck-integral-descent}) on integral descent of \'etale
maps. Finally, we prove~\ref{pro:equalizer-descent}, which provides the plan of the proof of
\ref{thm:witt-etale-general}. Aside from the language of descent, only these
three results will be used outside this section.

\spacesubsec{Language}

\subsection{} \emph{Fibered categories.}
Let $\cata$ be a category with fibered products. Let $\catb$ be a category fibered over $\cata$.
(See~\cite{TDTE-I}, A.1.1, or \cite{SGA1}, VI.6.1.) For any object $S$ of $\cata$, let $\catb_S$
denote the fiber of $E$ over $S$. Let us say that a map $q\:T\to S$ in $\cata$ is an
$\catb$-equivalence if $q^*\:\catb_S\to\catb_T$ is an equivalence of categories, and let us say
that $q$ is a \emph{universal $\catb$-equivalence} if for any map $S'\to S$ in $\cata$, the base
change $q'\:S'\times_S T\to S'$ is an $\catb$-equivalence.

For the applications in the next section, the reader can take
\begin{equation}
	\begin{split} \label{fibered-cat-example}
		\cata &= \text{the category of affine schemes,} \\
		\catb &= \text{the fibered category over $\cata$ where $\catb_S$ is the category of affine} \\
		 	& \phantom{=\ } 
				\text{\'etale $S$-schemes and the functors $q^*$ are given by base change.}
	\end{split}
\end{equation}
Then any closed immersion defined by a nil ideal
is a universal $\catb$-equivalence (EGA IV 18.1.2 \cite{EGA-no.32}).

\subsection{} \emph{Composition notation.}
Let $S$ be an object of $\cata$,
and let $\cata_{S\times S}$ denote the category of objects
over $S\times S$.  That is, an object of $\cata_{S\times S}$ is a pair $(T,\pp{T})$, where
$T$ is an object of $\cata$ and
$\pp{T}$ is a map $T\to S\times S$, called its structure map;
a morphism is a morphism in $\cata$ commuting with the maps
to $S\times S$. For such an object,
let $\pp{T,1},\pp{T,2}$ denote the composition of the structure map $T\to S\times S$ with
the projections $\pr_1,\pr_2\:S\times S\to S$. ($\pp{T,1}$ is the `source,' and $\pp{T,2}$ is
the `target'.)  We will often abusively leave $\pp{T}$
implicit and say that $T$ is an object of $\cata$.

Let $1_S$ denote the object $(S,\Delta)$ of $\cata_{S\times S}$, where $\Delta:S\to S\times S$ is
the diagonal map.

Given two objects $T,U\in\cata_{S\times S}$, define $TU\in\cata_{S\times S}$ 
as follows.  As an object of $\cata$, it is the fibered product
	\begin{equation}
		\xymatrix{
		TU \ar^{\pr_1}[r] \ar_{\pr_2}[d] & T \ar^{\pp{T,2}}[d] \\
		U \ar^{\pp{T,1}}[r] & S.
		}
	\end{equation}
We give $TU$ the structure of an object of $\cata_{S\times S}$ with the map
	\begin{equation}
		\xymatrix@C=85pt{
		TU = T\times_{S} U \ar^-{(\pp{T,1}\circ\pr_1,\pp{U,2}\circ\pr_2)}[r]
			& S\times S.
		}
	\end{equation}

\subsection{} \emph{Category objects and equivalence relations.}
A category object over $S$ is an object $R\in\cata_{S\times S}$ together with maps
	\begin{equation} \label{eq:groupoid-structure-maps}
		e_R\:1_S \to R, \quad c_R\:RR\to R
	\end{equation}
in $\cata_{S\times S}$ (called identity and composition)
satisfying the usual identity and associativity axioms in the definition
of a category.  

A morphism $f\:R\to R'$ of such category objects defined to be a morphism in $\cata_{S\times S}$
satisfying the functor axioms, that is, such that
	$$
	f\circ e_R = e_{R}\circ f \quad \text{and} \quad c_{R'}\circ ff = f\circ c_R,
	$$
where $ff$ denotes the map $RR\to R'R'$ induced by $f$.

A category-object structure on a subobject $R\subseteq S\times S$ is 
a property of $R$ in that when it exists, it is unique.
One might say that $R$ is a reflexive transitive relation on $S$.
We say $R$ is an equivalence relation on $S$ if, in addition,
the endomorphism $(\pr_2,\pr_1)$
of $S\times S$ that switches the two factors restricts
to a map 
	$$
	\inv\:R\to R
	$$
(which is of course unique when it exists).

\subsection{} \emph{Pre-actions (gluing data).}
\label{subsec:pre-actions}
Let $T$ be an object of $\cata_{S\times S}$.
A \emph{pre-action} of $T$ on an object $X\in\catb_{S}$ is defined to be an isomorphism
	\begin{equation}
		\varphi\:\pp{T,2}^*(X) \longisomap \pp{T,1}^*(X)
	\end{equation}
in $\catb_T$. A pre-action is also called a gluing datum on $X$ 
relative to the pair of maps $(\pp{T,1},\pp{T,2})$.
(Actually, Grothendieck \cite{TDTE-I}, A.1.4, calls $\varphi^{-1}$ the gluing datum.)
Let $\PreAct{T}{X}$ denote the set of pre-actions of $T$ on $X$.
Any map $T\to T'$ in $\cata_{S\times S}$ naturally induces a map
	$$
	\PreAct{T'}{X} \longmap \PreAct{T}{X}.
	$$

If $f\:X\to X'$ is a morphism in $\catb_S$ between objects $X,X'$ with pre-actions
$\varphi,\varphi'$, then we say $f$ is \emph{$T$-equivariant} if the following diagram
commutes:
$$
\xymatrix@C=90pt@R=50pt@!0{
\pp{T,2}^*(X) \ar^{\pp{T,2}^*(f)}[r]\ar^{\varphi}[d]
	& \pp{T,2}^*(X')\ar^{\varphi'}[d] \\
\pp{T,1}^*(X)\ar^{\pp{T,1}^*(f)}[r]
	& \pp{T,1}^*(X').
}
$$
In this way, the objects of $\catb_S$ equipped with a pre-action of $T$ form a category.

\subsection{} \emph{Actions.}
\label{subsec:actions}
Now let $R$ be a category object over $S$.  
An \emph{action} of $R$ on $X$ is defined
to be a pre-action $\varphi$ of $R$ on $X$ such that the diagram
	$$
	\xymatrix@C=50pt@R=50pt@!0{
	e^*\pp{R,2}^*(X) \ar^{e^*(\varphi)}[rr]\ar@{=}[dr]
		&
		& e^*\pp{R,1}^*(X) \ar@{=}[dl] \\
	& \id_S^*(X) \\
	}
	$$
and the diagram
	$$
	\xymatrix@C=90pt@R=50pt@!0{
		& c^*\pp{R,2}^*(X) \ar^{c^*(\varphi)}[r] \ar@{=}[dl]
		& c^*\pp{R,1}^*(X) \\
	\pr_2^*\pp{R,2}^*(X) \ar^{\pr_2^*(\varphi)}[dr]
		&
		&
		& \pr_1^*\pp{R,1}^*(X) \ar@{=}[ul] \\
		& \pr_2^*\pp{R,1}^*(X) \ar@{=}[r]
		& \pr_1^*\pp{R,2}^*(X), \ar^{\pr_1^*(\varphi)}[ur]
	}
	$$
commute.  Here, $\pr_1$ and $\pr_2$ denote the projections $RR\to R$ onto the first
and second factors, and the morphisms represented by equality signs are the isomorphisms
induced by the canonical
structure maps $(g\circ f)^*\isomap f^*\circ g^*$ (notated 
$c_{f,g}$ in \cite{TDTE-I}, A.1.1(ii))
of the fibered category $\catb$ corresponding to the equalities
	$$
	\pp{R,2}\circ e = \id_S = \pp{R,1}\circ e
	$$
and
	$$
	\pp{R,2}\circ c = \pp{R,2}\circ\pr_2,
	\quad \pp{R,1}\circ\pr_2 = \pp{R,2}\circ\pr_1, \quad \pp{R,1}\circ c = \pp{R,1}\circ\pr_1.
	$$
We will often use the following more succinct, if slightly abusive, expressions of
the commutativity of the diagrams above:
	\begin{equation} \label{eq:groupoid-action-def}
		e^*(\varphi) = \id_X, \quad c^*(\varphi) = (\pr_1^*\varphi)\circ(\pr_2^*\varphi).
	\end{equation}

Let $\Act{R}{X}$ denote the set of actions of $R$ on $X$.
A morphism $R\to R'$ of category objects induces a map
	$$
	\Act{R'}{X} \longmap \Act{R}{X}
	$$
in the obvious way.

Last, note that if $R$ is an equivalence relation, then the diagram 
	$$
\xymatrix@C=90pt@R=50pt@!0{
	\inv^*\pp{R,2}^*(X) \ar^{\inv^*(\varphi)}[r] \ar@{=}[d]
		& \inv^*\pp{R,1}^*(X) \ar@{=}[d]\\
	\pp{R,1}^*(X) \ar^{\varphi^{-1}}[r] 
		& \pp{R,2}^*. 
	}
	$$
commutes. This follows immediately from~(\ref{eq:groupoid-action-def}).
The abbreviated version is
	\begin{equation} \label{eq:groupoid-action-inverse}
		\inv^*(\varphi) = \varphi^{-1}.
	\end{equation}

\subsection{} \emph{Descent data.}
Let $q\:S'\to S$ be a map in $\cata$, and put
	$$
	R(S'/S)=S'\times_S S'.
	$$
View $R(S'/S)$ as an object in $\cata_{S'\times S'}$ by taking $\pi_{R(S'/S)}$
to be the evident monomorphism
	$$
	R(S'/S) = S'\times_S S' \longmap S'\times S'
	$$ 
Then $R(S'/S)$ is an equivalence relation on $S'$. An action $\varphi$ of
$R(S'/S)$ on an object $X'$ of $\catb_{S'}$ is also called a descent datum on $X'$ from $S'$ to 
$S$. (Again, it is actually $\varphi^{-1}$ that is called the descent datum in \cite{TDTE-I}.)
We might call $R(S'/S)$ the descent, or Galois, groupoid of the map $q\:S'\to S$.

Because the two compositions $R(S'/S)=S'\times_S S'\rightrightarrows S'\to S$ are equal,
for any object $X\in\catb_S$, the object $q^*(X)$ of $\catb_{S'}$ has a canonical pre-action
of $R(S'/S)$, and it is easy to check that this is an action.
We say that $q$ is a \emph{descent map} 
for the fibered category $\catb$ if the functor
from $\catb_S$ to the category of objects of $\catb_{S'}$ with an $R$-action is
fully faithful.  We say it is an \emph{effective descent map} if it is
an equivalence.

\subsection{} \emph{When gluing data is descent data.}
\label{subsec:when-gluing-is-descent}
Now suppose we have a diagram
	\begin{equation} \label{diag:gluing-descent}
		\displaycofork{S''}{S'}{S}
	\end{equation}
in $\cata$ such that the two compositions $S''\rightrightarrows S$ are equal. 
The universal property of products gives a map 
	$$
	S''\longmap S'\times_S S' = R(S'/S).
	$$
For any object $X'\in \catb_{S'}$, 
this map induces a function
	$$
	\Act{R(S'/S)}{X'} \longmap \PreAct{S''}{X'}.
	$$
Let us say that \emph{gluing data on $X'$ is descent data} relative to
the diagram~(\ref{diag:gluing-descent}) when this map is a bijection.

\spacesubsec{Gluing two objects}

Here we spell out in perhaps excessive detail some basic facts about equivalence relations
on disjoint unions which are $\catb$-trivial (though not necessarily trivial) on each factor.

From now on, let $\cata$ denote the category of affine schemes, schemes, or
algebraic spaces.  (We only need some weak hypotheses on coproducts in $\cata$, 
but let us not bother to determine which ones we need.)

\subsection{} \emph{Equivalence relations on a disjoint union.}
\label{subsec:gluing-two-objects}
Suppose $S$ is a coproduct $S_a + S_b$ of two objects $S_a,S_b\in\cata$.
(We use the symbols $a,b$ to index the summands only to emphasize their distinction
from the symbols $1,2$ that index the factors in the product $S\times S$.)
Let $R$ be an equivalence relation on $S$,
and let $R_{ij}$ denote $R\times_{S\times S}(S_i\times S_j)$, for any $i,j\in\{a,b\}$.
Let $\pp{R_{ij},1}$ denote the evident composition
	$$
	R_{ij} = R\times_{S\times S}(S_i\times S_j) \longlabelmap{\pr_1} S_i 		
	$$
and $\pp{R_{ij},2}$ the analogous map $R_{ij}\to S_j$.
We will sometimes view $R_{ij}$ as an object of $\cata_{S\times S}$
using the induced map $R_{ij} \to S_i\times S_j \to S\times S$.

Let $e_i\:S_i\to R_{ii}$ and $c_{ijk}\:R_{ij}R_{jk}\to R_{ik}$
and $\inv_{ij}\:R_{ij}\to R_{ji}$ denote the evident restrictions of $e$ and $c$ and $\inv$.

\subsection{} \emph{Actions over a disjoint union.}
For any object $X$ over $S$,
write $X_a= S_a\times_S X$ and $X_b= S_b\times_S X$.

For any pre-action 
	\begin{equation}
		\varphi\: \pp{R,2}^*X \longmap \pp{R,1}^*X,
	\end{equation}
of $R$ on $X$, let us write
$\varphi_{ij}$ for the restriction of $\varphi$ to $R_{ij}$.
In order for this pre-action to be an action, it is necessary 
and sufficient that for all $i,j,k\in\{a,b\}$ we have 
\begin{align}
		e_i^*(\varphi_{ii}) &= \id_{X_i} \text{ and }
		 	\label{eq:groupoid-relations-identity} \\
		c_{ijk}^*(\varphi_{ik}) &= \pr_1^*(\varphi_{ij})\circ\pr_2^*(\varphi_{jk}).
			\label{eq:groupoid-relations-associativity}
\end{align}
This is just a restatement of (\ref{eq:groupoid-action-def}), summand by summand.
In that case, \ref{eq:groupoid-action-inverse} becomes
	\begin{equation} \label{eq:groupoid-inverse}
		\inv_{ij}^*(\varphi_{ji}) = \varphi_{ij}^{-1}.
	\end{equation}

\begin{proposition}\label{pro:abstract-descent=gluing}
	Let $R$ be an equivalence relation on $S=S_a+S_b$ such that
	for $i=a,b$, the map $e_i\:S_i\to R_{ii}$ is a universal $\catb$-equivalence.
	Then for any object $X\in\catb_S$,
	the map 
		$$
		\xymatrix@C=35pt{
		\Act{R}{X} \ar^-{\varphi\mapsto\varphi_{ba}}[r]
		 	& \PreAct{R_{ba}}{X}
		}
		$$
	is a bijection.
\end{proposition}

\begin{proof} 
Let us first show injectivity.  
Let $\varphi$ and $\varphi'$ be actions of $R$ on $X$
such that $\varphi_{ba} = \varphi'_{ba}$.
We need to show that this implies $\varphi_{ij}=\varphi'_{ij}$ for all $i,j\in\{a,b\}$.
Consider each case separately.  For $ij=ba$, it is true by assumption.
When $ij=ab$, equation~(\ref{eq:groupoid-inverse})
and the given equality $\varphi_{ba}=\varphi'_{ba}$, imply
	$$
	\varphi_{ab}=\inv_{ba}^*(\varphi_{ab})^{-1} = \inv_{ba}^*(\varphi'_{ab})^{-1} = \varphi'_{ab}.
	$$
When $i=j$, since $e_i$ is
an $\catb$-equivalence, it is enough to show $e_i^*(\varphi_{ii})=e_i^*(\varphi'_{ii})$.
But by~(\ref{eq:groupoid-relations-identity}), we have
	$$
	e_i^*(\varphi_{ii}) = \id_{X_i} = e_i^*(\varphi'_{ii}).
	$$
Therefore $\varphi=\varphi'$, which proves injectivity.

Now consider surjectivity.  Let $\varphi_{ba}$ be a pre-action of $R_{ba}$ on $X$.
Define 
	\begin{equation}\label{eq:desc=glue-1}
		\varphi_{ab}=\inv_{ab}^*(\varphi_{ba})^{-1}
	\end{equation}
and for $i=a,b$ define
$\varphi_{ii}$ to be the map such that 
	\begin{equation}\label{eq:desc=glue-2}
		e_i^*(\varphi_{ii})=\id_{X_i},
	\end{equation}
which exists and is unique because $e_i$ is an $\catb$-equivalence.
We need to check that the pre-action
$\varphi=\varphi_{aa}+\varphi_{ab}+\varphi_{ba}+\varphi_{bb}$ of $R$ on $X$ is actually
an action.  To do this, we will verify the 
relations~(\ref{eq:groupoid-relations-identity}) and 
(\ref{eq:groupoid-relations-associativity}).

The identity axiom (\ref{eq:groupoid-relations-identity}) holds because it
is the defining property (\ref{eq:desc=glue-2}) of $\varphi_{ii}$.

Now consider the associativity axiom (\ref{eq:groupoid-relations-associativity}) 
for the various possibilities for $ijk$.  
Since $i,j,k\in \{a,b\}$, two of $i,j,k$ must be equal.

If $i=j$, the composition $f$ 
	$$
	\xymatrix{
	R_{jk} \ar^-{\pr_2^{-1}}_-{\sim}[r] & S_{jj}R_{jk} \ar^-{e_j\times \id}[r] & R_{jj}R_{jk}
	}
	$$  
is an $\catb$-equivalence, because it is a base change of
the universal $\catb$-equivalence $e_j$.  
Therefore it is enough to show 
	\begin{equation} \label{eq:abstract-gluing-1}
		f^*c_{jjk}^*(\varphi_{jk}) = f^*\pr_1^*(\varphi_{jj})\circ f^*\pr_2^*(\varphi_{jk}).
	\end{equation}
By the equality $\pr_1\circ f=e_j\circ \pp{R_{jk},1}$ and (\ref{eq:desc=glue-2}), we have 
	$$
	f^*\pr_1^*(\varphi_{jj}) = \pp{R_{jk},1}^*e_j^*(\varphi_{jj}) = \pp{R_{jk},1}^*(\id_{X_j})
		= \id.
	$$
On the other hand, by $c_{jjk}\circ f=\id_{R_{jk}}=\pr_2\circ f$, we have
$f^*c_{jjk}^*(\varphi_{jk}) = f^*\pr_2^*(\varphi_{jk})$.
Equation~(\ref{eq:abstract-gluing-1}) then follows.

The case $j=k$ is similar to the case $i=j$.  (Or apply $\inv$ to the case $i=j$.)

Last, suppose $i=k$. 
The following diagram is easily checked to be cartesian:
	\begin{equation} \label{diag:groupoid-cartesian-square}
		\xymatrix@C=90pt@R=50pt@!0{
		R_{ij} \ar^-{(\id_{R_{ij}},\inv_{ij})}[r]\ar_{\pp{R_{ij},1}}[d]
			& R_{ij}R_{ji} \ar^{c_{iji}}[d] \\
		S_i \ar^-{e_i}[r]
			& R_{ii}.
		}
	\end{equation}
(This is just another expression of the existence and uniqueness of inverses in a groupoid.)
Since $e_i$ is a universal $\catb$-equivalence, $(\id_{R_{ij}},\inv_{ij})$ 
is an $\catb$-equivalence.  So it is enough to show axiom 
(\ref{eq:groupoid-relations-associativity}) after applying
$(\id_{R_{ij}},\inv_{ij})^*$, that is, to show
	\begin{equation} \label{eq:abstract-gluing-2}
		(\id_{R_{ij}},\inv_{ij})^*c_{iji}^*(\varphi_{ii}) = 
			(\id_{R_{ij}},\inv_{ij})^*\pr_1^*(\varphi_{ij})
			\circ (\id_{R_{ij}},\inv_{ij})^*\pr_2^*(\varphi_{ji}).
	\end{equation}
By the commutativity 
of~(\ref{diag:groupoid-cartesian-square}) and (\ref{eq:desc=glue-2}),
we have 
	$$
	(\id_{R_{ij}},\inv_{ij})^*c_{iji}^*(\varphi_{ii})=\pp{R_{ij},1}^*e_i^*(\varphi_{ii})
		=\pp{R_{ij},1}^*(\id_{X_i})=\id.
	$$  
Combining this with the equation
$\varphi_{ji}=\inv_{ji}^*(\varphi_{ij})^{-1}$, equation~(\ref{eq:abstract-gluing-2})
reduces to
	$$
	(\id_{R_{ij}},\inv_{ij})^*\pr_1^*(\varphi_{ij})
		=(\id_{R_{ij}},\inv_{ij})^*\pr_2^*\inv_{ji}^*(\varphi_{ij}).
	$$  
But this holds because we have
	\begin{equation*}
		\pr_1\circ (\id_{R_{ij}},\inv_{ij}) = \id_{R_{ij}} 
			=\inv_{ji}\circ \inv_{ij}=\inv_{ji}\circ \pr_2\circ (\id_{R_{ij}},\inv_{ij}).
	\end{equation*}
	
Therefore the equations in~(\ref{eq:groupoid-relations-associativity}) hold for all $i,j,k$, and
so the pre-action is an action.
\end{proof}

\spacesubsec{Grothendieck's theorem}

Recall that a map $\Spec B\to \Spec A$ of affine schemes is said to be integral if the
corresponding ring map $A\to B$ is integral (and not necessarily injective).

\begin{theorem}\label{thm:Grothendieck-integral-descent}
	Every surjective integral map $Y\to X$ of affine schemes
	is an effective descent map for the fibered category $\catb$ over $\cata$
	of~(\ref{fibered-cat-example}).
\end{theorem}

This theorem is proven in SGA 1 IX 4.7~\cite{SGA1} up to two details. First, the argument given
there covers only
morphisms $Y\to X$ which are finite and of finite presentation; and second, the
statement there has no affineness in the assumptions or in the conclusion. The first point can be
handled by a standard limiting argument (or one can apply Theorem 5.17 plus Remark 2.5(1b) in 
Rydh~\cite{Rydh:Submersions}). The second point can be handled with Chevalley's theorem; the form
most convenient here would the final one, Rydh's~\cite{Rydh:Noetherian-approximation} Theorem
(8.1), which is free of noetherianness, separatedness, finiteness, and scheme-theoretic
assumptions.

\spacesubsec{Gluing and descent of \'etale algebras}

\begin{proposition}\label{pro:equalizer-descent}
Consider a diagram of rings
	\begin{equation} \label{diag:equalizer-descent}
	\xymatrix{
	B \ar^{d}[r]
		& B'\rightlabelxyarrows{\ha}{\hb} 		
		& B'' \\
	A \ar_{\ff}[r]\ar^e[u]
		& A' \rightlabelxyarrows{\ga}{\gb}
			\ar^{e'}[u]
		& A''\ar_{e''}[u]
	}
	\end{equation}
such that 
$ h_i\circ e'=e''\circ g_i$, for $i=1,2$.
Also assume the following properties are satisfied:
	\begin{enumerate}
		\item the two parallel right-hand squares are cocartesian,
		\item both rows are equalizer diagrams,
		\item relative to the lower row, gluing data on any \'etale $A'$-algebra
			is descent data,
		\item $f$ satisfies effective descent for the fibered category of \'etale algebras,
			and
		\item $e'$ is \'etale.
	\end{enumerate}
Then $e$ is \'etale and the left-hand square is cocartesian.
\end{proposition}

Note that when we use the language of descent in the category of rings (as in (c)
and (d)), we understand that it refers to the corresponding statements in the opposite category.

\begin{proof}
	Property (a) equips the \'etale $A'$-algebra $B'$ 
	with gluing data $\varphi$ relative to $(\ga,\gb)$.  Indeed, take $\varphi$ to be the
	composition
		$$
		A''\tn_{g_1,A'} B' \longisomap B'' \longisomap A''\tn_{g_2,A'} B'.
		$$
	By property (c), this gluing data comes from unique descent data relative
	to $f$.
	Therefore by (d) and (e), the $A'$-algebra
	$B'$ descends to an \'etale $A$-algebra $C$.
	
	Now apply the functor $C\tn_A\vbl$ to the lower row of
	diagram~(\ref{diag:equalizer-descent}).  By (a) and the definition
	of descent, the result can be identified with the sequence
		\[
		\xymatrix{
		C \ar[r]
			& B' \rightlabelxyarrows{\ha}{\hb} 		
			& B''. \\
		}
		\]
	This sequence is also an equalizer diagram, because the lower row
	of~(\ref{diag:equalizer-descent}) is an equalizer diagram, by (b), and because
	$C$ is \'etale over $A$ and hence flat.  Again by (b), the upper
	row of~(\ref{diag:equalizer-descent}) is an equalizer diagram, and
	so we have $C=B$.  Therefore, $B$ is an \'etale $A$-algebra and the left-hand
	square is cocartesian.
\end{proof}

\section{Ghost descent in the single-prime case}
\label{sec:affine-ghost-descent-single-prime-case}

We return to the notation of~\ref{subsec:affine-general-notation}.
Suppose $\ptst$ consists of a single maximal ideal $\m$, and fix an integer $n\geq 1$.
Write $W_n=W_{R,\m,n}$, and so on.
Let $A$ be an $R$-algebra, and let $\alpha_n$ denote the map
		\begin{equation}
			W_n(A) \xright{\alpha_n} W_{n-1}(A) \times A
		\end{equation}
given by the canonical projection on the factor $W_{n-1}(A)$ and the $n$-th 
ghost component $\gh{n}$ on the factor $A$.
Let $I_n(A)$ denote the kernel of $\alpha_n$.
For example, if $\m$ is generated by $\pi$, then in terms of the Witt components, we have
	\begin{equation} \label{eq:concrete-description-of-kernel-alpha}
		I_n(A) = \setof{(0,\dots,0,a)_{\pi}\in A^{[0,n]}}{\pi^n a = 0}.
	\end{equation}

\begin{proposition}\label{pro:alpha-map} 
	We have the following:
	\begin{enumerate}
		\item \label{part:alpha-map-a}
			$\alpha_n$ is an integral ring homomorphism.
		\item \label{part:alpha-map-b}
			The kernel $I_n(A)$ of $\alpha_n$ is a square-zero ideal.
		\item \label{part:alpha-map-c}
			If $A$ is $\m$-flat, then $\alpha_n$ is injective. 
		\item \label{part:alpha-map-d}
			The diagram
				$$
				\xymatrix{
				W_{n-1}(A) \ar^-{\rgh{n}}[r]
					& A/\m^n A \\
				W_n(A) \ar^-{\gh{n}}[r]\ar@{->>}[u]
					& A, \ar@{->>}[u]
				}
				$$
			where the vertical maps are the canonical ones, is cocartesian.
		\item \label{part:alpha-map-e}
			View $A$ as a $W_n(A)$-algebra by the map $\gh{n}\:W_n(A)\to A$.
			Then every element of the kernel of the multiplication map
			 	$$
				A\tn_{W_n(A)} A \longmap A
				$$
			is nilpotent.
		\item \label{part:alpha-map-f}
			In the diagram
				\begin{equation} \label{diag:alpha-equalizer}
					\displaylabelfork{W_n(A)}{\alpha_n}{W_{n-1}(A)\times A}
						{\rgh{n}\circ\pr_1}{\overline{\pr}_2}{A/\m^n A,}
				\end{equation}
			where $\overline{\pr}_2$ denotes the reduction of $\pr_2$ modulo $\m^n$, 
			the image of $\alpha_n$ agrees with the equalizer of 
				${\rgh{n}\circ\pr_1}$ and ${\overline{\pr}_2}$.
	\end{enumerate}
\end{proposition}

\begin{proof}
(\ref{part:alpha-map-a}): 
It is enough to show that each factor of $W_{n-1}(A)\times A$ is integral 
over $W_{n}(A)$.  The first factor is a quotient ring, and hence integral.  Now consider 
an element $a\in A$.  Then $a^{q^{n}}$
is the image in $A$ of the Teichm\"uller lift $[a]\in W_{n}(A)$. 
(See~\ref{subsec:teich-1}.)  Therefore the 
second factor is also integral over $W_{n}(A)$.

(\ref{part:alpha-map-b}): 
It suffices to show this 
after base change to $R[1/\m]\times R_{\m}$.
Therefore, by~\ref{pro:Zariski-localization-on-R}, we may assume $\m$ is generated
by a single element $\pi$.  Then an element of the kernel
of $\alpha_n$ will be of the form $V_\pi^n[a]=(0,\dots,0,a)_{\pi}$, where $\pi^na=0$.
On the other hand, by~(\ref{eq:Verschiebung-product-new}) we have
	$$
	(V_\pi^n[a])(V_\pi^n[b]) = \pi^nV_\pi^n[ab] = (0,\dots,0,\pi^n ab)_{\pi} = 0.
	$$

(\ref{part:alpha-map-c}): We have $(\gh{\leq n-1}\times \id_A)\circ\alpha_n=\gh{\leq n}$.
Since $A$ is $\m$-flat, the map $\gh{\leq n}$ is injective (\ref{pro:ghost-map-inj}), 
and hence so is $\alpha_n$.

(\ref{part:alpha-map-d}):
As above, it is enough by~\ref{pro:Zariski-localization-on-R} to assume $\m$ is
generated by a single element $\pi$.  Then we have
	$$
	A\tn_{W_n(A)} W_{n-1}(A) = A\tn_{W_n(A)}W_n(A)/V^nW_n(A) = A/\gh{n}(V^nW_n(A))A.
	$$
Examining the Witt polynomials~(\ref{eq:witt-comp}) shows $\gh{n}(V^nW_n(A))=\pi^n A$.

(\ref{part:alpha-map-e}): 
Again, by~\ref{pro:Zariski-localization-on-R} 
we may assume $\m$ is generated by a single element $\pi$.
To show every element $x\in I$ is nilpotent, it is enough to restrict $x$
to a set of generators.  Therefore it is enough to show $(1\tn a-a\tn 1)^{q^n}=0$ for
every element $a\in A$.  

Now suppose that, for $j=0,\dots,q^n$, we could show
	\begin{equation} \label{eq:alpha-map-6b}
		\binom{q^n}{j}a^j \in \im(\gh{n}).
	\end{equation}
Then we would have
	\begin{multline*}
		(1\tn a-a\tn 1)^{q^n} = \sum_{j} (-1)^j\binom{q^n}{j} a^j\tn a^{q^n-j} \\
			= \sum_j (-1)^j\tn \binom{q^n}{j} a^j a^{q^n-j} 
			= (1\tn a- 1\tn a)^{q^n}=0,
	\end{multline*}
which would complete the proof.  So let us show~(\ref{eq:alpha-map-6b}).
	
Let $f=\ord_p(q)$ and $i=\ord_p(j)$. 
Then we have 
	$$
	\ord_p\binom{q^n}{j} = \ord_p(q^nj^{-1}) + \ord_p\binom{q^n-1}{j-1}
			\geq nf-i.
	$$
It follows that $\binom{q^n}{j}a^j$ is an $R$-linear multiple of $\pi^{nf-i}a^j$. 
Since $\gh{n}$ is an $R$-algebra map, it is therefore enough to show
	\begin{equation} \label{eq:alpha-map-6d}
	 	\pi^{nf-i}a^j \in \im(\gh{n}).  
	\end{equation}

Now, for $b\in A$ and $s=0,\dots,n$, we have $\pi^{n-s} b^{q^{s}}=\gh{n}(V^{n-s}_\pi[b])$,
and therefore $\pi^{n-s} b^{q^{s}}$ is in the image of $\gh{n}$.
So to show~(\ref{eq:alpha-map-6d}), it is enough to 
find an integer $s$ and an element $b\in A$ such that 
$\pi^{n-s} b^{q^{s}}$ is an $R$-linear divisor of $\pi^{nf-i}a^j$.
In particular, it is sufficient for $b$ and $s$ to satisfy $b^{q^s}=a^j$ and $n-s\leq nf-i$.

Take $s$ to be the greatest integer at most $if^{-1}$. Then we have $q^s\mid j$; so if we set
$b=a^{j/q^s}\in A$, we have $b^{q^s}=a^j$. It remains to show $n-s\leq nf-i$. This is equivalent to
$n-if^{-1} \leq nf-i$, which is in turn equivalent to $(1-f)(n-if^{-1})\leq 0$. And this holds
because $1-f\leq 0$ and $n-if^{-1}\geq 0$. (Recall that $j\leq q^n$.) This completes the proof of
(e).

(\ref{part:alpha-map-f}): 
As above, we may assume that $\m$ can be generated by a single element $\pi$.  
For any element $a=(a_0,\dots,a_n)_{\pi}\in W_{n}(A)$, we have
	$$
	\alpha_n(a) = \Big((a_0,\dots,a_{n-1}),
		a_0^{q^n}+\cdots + \pi^{n-1} a_{n-1}^q +\pi^{n} a_{n}\Big).
	$$
Therefore an element $\big((a_0,\dots,a_{n-1}),b\big)\in W_{n-1}(A)\times A$
lies in the image of $\alpha_n$ if and only if
	$$
	a_0^{q^n}+\cdots + \pi^{n-1} a_{n-1}^q \equiv b \mod\m^n A,
	$$
which is exactly what we needed to show.
\end{proof}

\begin{corollary}\label{cor:gh-integral-and-surj-single-prime}
	For any $R$-algebra $A$, the ghost map
		$$
		\gh{\leq n}\:W_n(A) \longmap A^{[0,n]}
		$$
	is integral, and its kernel $J$ satisfies $J^{2^n}=0$.
\end{corollary}

\begin{proof}
	By \ref{pro:alpha-map} and induction on $n$.
\end{proof}

\begin{theorem}\label{thm:alpha-gluing=descent}
	\begin{enumerate}
		\item The map $\alpha_n$ is an effective descent map for the fibered category
			of \'etale algebras.
		\item Relative to the diagram
			\begin{equation} \label{diag:main-gluing-flat}
				\displaylabelfork{W_n(A)}{\alpha_n}{W_{n-1}(A)\times A}
					{\rgh{n}\circ\pr_1}{\overline{\pr}_2}{A/\m^nA,}			
			\end{equation}
			gluing data on any \'etale $W_{n-1}(A)\times A$-algebra is descent 
			data~(\ref{subsec:when-gluing-is-descent}).
		\item If $A$ is $\m$-flat, then for 
			any $A'$-algebra $B'$ equipped with gluing data $\varphi$,
			the descended $A$-algebra is the subring $B$ of $B'$ on which the following diagram 
			commutes:
				$$
				\xymatrix@R=20pt{
					& A/\m^nA\tn_{\rgh{n}\circ\pr_1}B' \\
				B' \ar^-{1\tn\id_{B'}}[ur]\ar_-{1\tn\id_{B'}}[dr] \\
					& A/\m^nA\tn_{\overline{\pr}_2}B'. \ar^{\varphi}[uu]
				}
				$$
	\end{enumerate}
\end{theorem}
\begin{proof}
	(a): This follows from 
	Grothendieck's theorem~\ref{thm:Grothendieck-integral-descent} 
	and~\ref{pro:alpha-map}(\ref{part:alpha-map-a})--(\ref{part:alpha-map-b}).

	(b): We will use~\ref{pro:abstract-descent=gluing},
	where $\cata$ and $\catb$ are as in (\ref{fibered-cat-example}).
	In the notation of~\ref{subsec:gluing-two-objects}, 
	put $S_a=\Spec W_n(A)$ and
	$S_b=\Spec A$.  Let $\Gamma$ be the equivalence relation $S\times_{\Spec W_n(A)} S$
	on $S$.  By~\ref{pro:alpha-map}(\ref{part:alpha-map-d}),
	we have $\Gamma_{ba}=\Spec A/\m^nA$.
	The map $e_a$ is an isomorphism because $W_{n-1}(A)$ is a quotient
	ring of $W_n(A)$.  The map $e_b$ is a nil immersion, 
	by~\ref{pro:alpha-map}(\ref{part:alpha-map-e}), and hence is an
	$\catb$-equivalence.
	Thus we can apply~\ref{pro:abstract-descent=gluing}, which says that
	a $\Gamma$-action is the same as a $\Gamma_{ba}$ pre-action.
	In other words, gluing data is descent data.
	
	(c): This will follow from~\ref{pro:equalizer-descent} once we verify the hypotheses.
	\ref{pro:equalizer-descent}(a)--(b) are clear; \ref{pro:equalizer-descent}(c)
	follows from (b) above; \ref{pro:equalizer-descent}(d) follows from (a) above;
	and	\ref{pro:equalizer-descent}(e) follows from the definition of $B$, for the top row
	of (\ref{diag:equalizer-descent}), and
	from~\ref{pro:alpha-map}(\ref{part:alpha-map-c}) and (\ref{part:alpha-map-f}),
	for the bottom row.
\end{proof}

\subsection{} \emph{Remark.}
\label{rmk:etale-alg-product-of-categories-description}
For any ring $C$, let $\EtAlg_C$ denote the category of \'etale $C$-algebras.
Then another way of expressing part (b) of this theorem is that the induced functor
	$$
	\EtAlg_{W_n(A)} \longmap  
		\EtAlg_{W_{n-1}(A)}\times_{\EtAlg_{A/\m^nA}}\EtAlg_{A}
	$$
is an equivalence. (Of course, the fibered product of categories is taken in
the weak sense.)
In particular, we can prove things about \'etale $W_n(A)$-algebras by induction on $n$.
This is the main technique in the proof of~\ref{thm:witt-etale-general}.
But it also seems interesting in its own right and will probably have applications
beyond the present paper.

\subsection{} \emph{Remark.}
\label{rmk:alpha-gluing=descent}
If we let $\wpre_n(A)$ denote the image of $\alpha_n$, then the induced diagram
		$$
		\displaylabelfork{\wpre_n(A)}{}{W_{n-1}(A)\times A}
			{\rgh{n}\circ\pr_1}{\overline{\pr}_2}{A/\m^nA.}
		$$
satisfies all the conclusions of the theorem above, regardless of whether $A$ is
$\m$-flat.

Indeed, it is an equalizer diagram by~\ref{pro:alpha-map}(\ref{part:alpha-map-f}) and the
definition of $\wpre_n(A)$; it is an effective descent map by~\ref{thm:alpha-gluing=descent}
and~\ref{pro:alpha-map}(\ref{part:alpha-map-b}). Last, because $\wpre_n(A)$ is the image of
$\alpha_n$, gluing (resp.\ descent) data relative to $\wpre_n(A)$ agrees with gluing (resp.\
descent) data relative to $W_n(A)$. In particular, gluing data relative to $\wpre_n(A)$ is descent
data.

\section{$W$ and \'etale morphisms}
\label{sec:W-and-etale-morphisms}

We return to the general context of~\ref{subsec:affine-general-notation}. In particular,
$\ptst$ is no longer required to consist of one ideal. 

\begin{lemma}\label{lem:etale-radicial-cartesian}
 	Consider a commutative square of affine schemes (or any schemes)
		$$
		\xymatrix{
		X \ar_{f}[d]		& X' \ar^{f'}[d]\ar_{g}[l] \\
		Y 					& Y', \ar_{h}[l]
		}
		$$
	and let $U$ be an open subscheme of $Y$. Suppose the following hold
		\begin{enumerate}
			\item $f$ and $f'$ are \'etale,
			\item the square above becomes cartesian after the base change $U\times_Y\vbl$,
			\item $g$ and $h$ become surjective and universally injective 
				after the base change $(Y-U)\times_Y \vbl$. 
		\end{enumerate}
	Then the square above is cartesian.
\end{lemma}
\begin{proof}
	Let $e$ denote the induced map $(g,f')\:X'\to X\times_Y Y'$.
	It is enough to show $e$ is \'etale, surjective,
	and universally injective (EGA IV (17.9.1)~\cite{EGA-no.32}).
	The composition of $e$ 	with $\pr_2\:X\times_Y Y'\to Y'$ is $f'$.  Because $f$ is 
	\'etale, so is its base change $\pr_2$.	Combining this with the \'etaleness of $f'$ implies
	that $e$ is \'etale (EGA IV 17.3.4~\cite{EGA-no.32}).

	The surjectivity and universal injectivity of $e$ can be checked after base change over $Y$
	to $U$ and to $Y-U$.  By assumption $e$ becomes an isomorphism after base
	change to $U$.  In particular, it becomes surjective and universally injective.
	
	Let $\bar{e},\bar{g},\bar{h}$ denote the maps $e,g,h$ pulled back from $Y$ to $Y-U$.
	Let $\bar{h}'$ denote the base change of $\bar{h}$ from $Y$ to $X$.
	Then, as above, we have $\bar{g}=\bar{h}'\circ\bar{e}$.
 	Since $\bar{h}$ is universally injective, so is $\bar{h}'$.
 	Combining this with the fact that $\bar{g}$ is universally injective, implies
	that $\bar{e}$ is as well (EGA I 3.5.6--7~\cite{EGA-no.4}). 
	Finally $\bar{e}$ is surjective since $\bar{h}'$ is injective and $\bar{g}$ is surjective. 
\end{proof}

\begin{theorem}\label{thm:witt-etale-general}
	For any \'etale map $\varphi\:A\to B$ and any element $n\in\nset$,
	the induced map $W_{R,\ptst,n}(\varphi)\:W_{R,\ptst,n}(A)\to W_{R,\ptst,n}(B)$ is \'etale.
\end{theorem}

\begin{proof}
By~\ref{cor:witt-iterate-generalized}, it is enough to assume $\ptst$ 
consists of a single maximal ideal $\m$. 
Also, it will simplify notation if we assume $\m$ is principal, generated by an element $\pi$.
We may do this by by~\ref{pro:Zariski-localization-on-R} and because it is enough to show 
\'etaleness after applying $R_{\m}\tn_R\vbl$ and $R[1/\m]\tn_R\vbl$. Let us write
$W_n=W_{R,\ptst,n}$.

We will reason by induction on $n$, the case $n=0$ being clear because $W_0$ is the identity
functor. So from now on, assume $n\geq 1$.

Let $\wpre_n(A)$ denote the image of $\alpha_n\:W_n(A)\to W_{n-1}(A)\times A$,
and let $\alphapre_n$ denote the induced injection $\wpre_n(A)\to W_{n-1}(A)\times A$.
Define $\wpre_n(B)$ and $\alphapre_n$ for $B$ similarly.

\emph{\underline{Step 1}}: $\wpre_n(B)$ is \'etale over $\wpre_n(A).$

	To do this, it suffices to verify conditions (a)--(e) of~\ref{pro:equalizer-descent}
	for the following diagram
		$$
		\xymatrix{
			\wpre_n(B) \ar^-{\alphapre_n}[r]
				& W_{n-1}(B) \times B\rightlabelxyarrows{\rgh{n}\circ\pr_1}{\overline{\pr}_2} 		
				& B/\m^n B \\
			\wpre_n(A) \ar^-{\alphapre_n}[r]\ar[u]
				& W_{n-1}(A) \times A \rightlabelxyarrows{\rgh{n}\circ\pr_1}{\overline{\pr}_2}
					\ar[u]
				& A/\m^n A,\ar[u]
		}
		$$
	where the vertical maps are induced by $\varphi$ and functoriality.
	We know~\ref{pro:equalizer-descent}(a) holds by induction.
	Conditions~\ref{pro:equalizer-descent}(c)--(d) hold by~\ref{thm:alpha-gluing=descent}
	(or \ref{rmk:alpha-gluing=descent}).
	Condition~\ref{pro:equalizer-descent}(e) was shown already 
	in~\ref{pro:alpha-map}(\ref{part:alpha-map-f}).
	Now consider~\ref{pro:equalizer-descent}(b).  It is clear that
	the square of $\overline{\pr}_2$ maps is cocartesian.  So, all that remains
	is to check that the square of $\rgh{n}\circ\pr_1$ maps is cocartesian.
	By induction, $W_{n-1}(B)$ is \'etale over $W_{n-1}(A)$, and so this follows
	from \ref{lem:etale-radicial-cartesian}, which we can apply
	by \ref{pro:Zariski-localization-on-R} and \ref{pro:W-nilpotent-at-p}. 
	
\emph{\underline{Step 2}}: $W_n(B)$ is \'etale over $W_n(A)$.

	By \ref{pro:alpha-map}(\ref{part:alpha-map-b}), the kernel $I_n(A)$ of the map 
	$\alpha_n\:W_n(A)\to\wpre_n(A)$ has square zero.  Therefore by EGA IV 18.1.2 \cite{EGA-no.32},
	there is an \'etale
	$W_n(A)$-algebra $C$ and an isomorphism $f\:C\tn_{W_n(A)}\wpre_n(A)\to \wpre_n(B)$.
	Now consider the square
		$$
		\xymatrix{
		C \ar[r]\ar^{d}@{-->}[dr]	& \wpre_n(B) \\
		W_n(A) \ar[u]\ar[r]			& W_n(B), \ar@{>>}[u]
		}
		$$
	where the upper map is the one induced by $f$
	and where $d$ will soon be defined.
	By~\ref{pro:alpha-map}(\ref{part:alpha-map-b}), 
	the kernel $I_n(B)$ of the right-hand map has square zero.  Therefore since $C$ is
	\'etale over $W_n(A)$, there exists a unique map $d$ making the diagram commute.
 	Let us now show that $d$ is an isomorphism.

	Because $C$ is \'etale and hence flat over $W_n(A)$,
	we have a commutative diagram with exact rows:
		$$
		\xymatrix{
		0 \ar[r]	& I_n(B) \ar[r]		& W_n(B)\ar[r] 	& \wpre_n(B)\ar[r]		& 0 \\
		0 \ar[r] 	& C\tn_{W_n(A)} I_n(A)\ar[r]\ar^{e}[u]			& C \ar[r]\ar^{d}[u] 
			& C\tn_{W_n(A)}\wpre_n(A) \ar[r]\ar_{\sim}^{f}[u]
			& 0.
		}
		$$
	So to show $d$ is an isomorphism, it is enough to show $e$ is an isomorphism.
	Because $I_n(A)$ is a square-zero ideal, the action of $W_n(A)$ on it factors through
	$\wpre_n(A)$.	Therefore, $e$ factors as follows:
		\begin{multline*}
		C\tn_{W_n(A)} I_n(A) = C\tn_{W_n(A)} \wpre_n(A) \tn_{\wpre_n(A)} I_n(A) \\
		\longlabelmap{f\tn\id}\wpre_n(B)\tn_{\wpre_n(A)} I_n(A) \longlabelmap{\xx} I_n(B),
		\end{multline*}
	Since $f$ is an isomorphism, it is enough to show $\xx$ is an isomorphism.
	
	Using the description (\ref{eq:concrete-description-of-kernel-alpha}) of $I_n$,
	the map $\xx$ can be extended to the following commutative diagram with 
	exact rows:
		$$
		\xymatrix{
		0 \ar[r]
			& I_n(B) \ar[r]
			& B \ar^{\cdot\pi^n}[r]
			& \pi^n B \ar[r]
			& 0\\
		0 \ar[r]
			& \wpre_n(B)\tn I_n(A) \ar[r]\ar^{\xx}[u]
			& \wpre_n(B)\tn A \ar^-{\cdot\pi^n}[r]\ar[u]^{(\pr_2\circ\alphapre_n)\cdot\varphi}
			& \wpre_n(B)\tn \pi^n A, \ar_{(\pr_2\circ\alphapre_n)\cdot\varphi}[u] \ar[r]
			& 0
		}
		$$
	where $\tn$ denotes $\tn_{\wpre_n(A)}$, for short.
	Therefore it is enough to show the right two vertical maps are isomorphisms,
	and to do this, it is enough to show the right-hand square in the diagram 
		\[
		\xymatrix@C=70pt@R=50pt@!0{
		W_n(B) \ar@{->>}[r]			& \wpre_n(B) \ar^-{{\pr}_2\circ\alphapre_n}[r]		 
			& B \\
		W_n(A) \ar@{->>}[r]\ar[u]	& \wpre_n(A) \ar[u]\ar^-{{\pr}_2\circ\alphapre_n}[r] 
			& A \ar[u]
		}
		\]
	is cocartesian.  We will do this by applying \ref{lem:etale-radicial-cartesian}, with
	$U=\Spec R[1/\m]\tn_R \wpre_n(A)$. 
	
	By step 1, condition \ref{lem:etale-radicial-cartesian}(a) holds.
	Now consider conditions \ref{lem:etale-radicial-cartesian}(b)--(c).
	By \ref{thm:alpha-gluing=descent}(\ref{part:alpha-map-b}), 
	the horizontal maps in the left-hand square
	have square-zero kernel.
	In particular, the scheme maps they induce are universal homeomorphisms.
	And by \ref{pro:Zariski-localization-on-R}, 
	they become isomorphisms after applying $R[1/\m]\tn_R\vbl$.
	Therefore it is enough to show \ref{lem:etale-radicial-cartesian}(b)--(c) hold
	for the perimeter of the diagram above. In this case,  
	\ref{lem:etale-radicial-cartesian}(b) follows from \ref{pro:Zariski-localization-on-R}, and
	\ref{lem:etale-radicial-cartesian}(c) follows from \ref{pro:W-nilpotent-at-p}.
\end{proof}

\subsection{} \emph{Remark.}
Observe that when $A$ is $\ptst$-flat, the proof terminates after step 1, which
is just an application of \ref{pro:equalizer-descent}.
Thus in the central case, the argument is not much more than \ref{pro:alpha-map} and some
general descent theory. 

\begin{corollary}\label{cor:affine-W-etale-base-change}
	Let $B$ an \'etale $A$-algebra, and let $C$ be any $A$-algebra.  Then
	for any $n\in\nset$, the induced diagram
		$$
		\xymatrix{
		W_{R,\ptst,n}(B) \ar[r]				& W_{R,\ptst,n}(B\tn_A C) \\
		W_{R,\ptst,n}(A) \ar[r]\ar[u]		& W_{R,\ptst,n}(C) \ar[u]
		}
		$$
	is cocartesian.
\end{corollary}
\begin{proof}
	By~\ref{cor:witt-iterate-generalized}, we can assume $\ptst$ consists of 
	a single ideal $\m$.
	The proof will be completed by \ref{lem:etale-radicial-cartesian}, once we check
	its hypotheses are satisfied.
	Condition (a) of \ref{lem:etale-radicial-cartesian} holds by~\ref{thm:witt-etale-general},
	condition (b) holds by \ref{pro:Zariski-localization-on-R} and \ref{pro:ghost-map-inj}, and 
	condition (c) holds by \ref{pro:W-nilpotent-at-p}.
\end{proof}

\subsection{} {\em $W_n$ does not generally commute with coproducts.}
\label{subsec:W-affine-does-not-commute-with-coproducts}
Almost anything is an example.
For instance, with the  $p$-typical Witt vectors,
$W_1(A\tn_\bZ A)$ is not isomorphic to $W_1(A)\tn_{W_1(\bZ)}W_1(A)$,
when $A$ is $\bF_p[x]$ or $\bZ[x]$.

\subsection{} \emph{$W$ does not generally preserve \'etale maps.}
Let $W$ denote $p$-typical Witt functor (non-truncated), and 
let $\varphi$ denote the evident inclusion $\bQ[x]\to\bQ[x^{\pm 1}]$, which is \'etale.
While $W(\varphi)$ is best viewed as a map of pro-rings, it is possible to view it as
a map of ordinary rings, and ask whether it is \'etale. It is not:
$W(\varphi)$ can be identified with 
$\varphi^{\bN}\:\bQ[x]^{\bN}\to\bQ[x^{\pm 1}]^{\bN}$, which is not \'etale because 
$\bQ[x^{\pm 1}]^{\bN}$ is not finitely generated as an $\bQ[x]^{\bN}$-algebra.
This is an elementary exercise.

\subsection{} \emph{Other truncation sets for the big Witt vectors.}
Some writers have considered more general systems of truncations for the big Witt functor
(\ref{subsec:agreement-with-classical-witt-vectors}). See
Hesselholt--Madsen~\cite{Hesselholt-Madsen:cyclic}, subsection 4.1, for example. Given a finite set
$T$ of positive integers closed under extraction of divisors, they define an
endofunctor $W_T$ of the category of rings. When $T$ consists of all the
divisors of some integer $d\geq 1$, then $W_T$ agrees with our $W_{\bZ,\ptst,n}$, where $\ptst$
consists of the maximal ideals $\m\subset\bZ$
that contain $d$ and where $n_{\m}=\ord_{\m}(d)$. Thus
the two systems of truncations are cofinal with respect to each other.

The functors $W_T$ also preserve \'etale maps. Indeed, it is enough to show that the base change to
$\bZ[1/T]$ and to $\bZ_{(p)}$, for each prime $p\in T$, is \'etale. (See EGA IV (17.7.2)(ii).)
Applying the identity $W_T(A)[1/p]=W_T(A[1/p])$, which can be established by looking at the graded
pieces of the Verschiebung filtration, it is enough to consider $\bZ[1/T]$-algebras and
$\bZ_{(p)}$-algebras. In the either case, $W_T(A)$ is simply a product of $p$-typical Witt rings
$W_n(A)$ for various primes $p$ and lengths $n$ (see \cite{Hesselholt-Madsen:cyclic}, (4.1.10)), in
which case the result follows from~\ref{thm:witt-etale-general}, or van der Kallen's original
theorem~\cite{van-der-Kallen:Descent}, (2.4).

\bibliography{references}
\bibliographystyle{plain}

\end{document}